\documentclass[12pt]{article}
\usepackage{amsmath,amssymb,amscd}
\newenvironment{demo}[1]{\par\smallskip\par\begin{trivlist}
\item[]{\bf #1}\ }{\end{trivlist}\par\smallskip\par}
\newcommand{\Proof}{\begin{demo}{{\it Proof.\ }}}
\newcommand{\qed}{\end{demo}}
\newcommand{\toy}{\ \rule[0em]{0.5ex}{1.8ex}}
\newcommand{\QED}{\toy\end{demo}}
\newtheorem{thm}{Theorem}[section]
\newtheorem{prop}[thm]{Proposition} 
\newtheorem{cor}[thm]{Corollary}
\newtheorem{lem}[thm]{Lemma}
\newtheorem{remark}[thm]{Remark}

\makeatletter \@addtoreset{equation}{section} \makeatother
\makeatletter
\ifnum\@ptsize=0

\else

\fi
\makeatother

\renewcommand{\a}{\alpha}
\renewcommand{\b}{\beta}

\newcommand{\dl}{\delta}

\newcommand{\lm}{\lambda}

\newcommand{\ve}{\varepsilon}
\newcommand{\D}{\Delta}

\newcommand{\nn}{\nonumber}

\begin{document}
\title{
Large deviation principle for certain spatially lifted Gaussian rough path
}
\author{Yuzuru INAHAMA 
\footnote{
Graduate School of Mathematics,   Nagoya University
%\\
Furocho, Chikusa-ku, Nagoya 464-8602, JAPAN.
%\\
E-mail:~\tt{inahama@math.nagoya-u.ac.jp}
}
}
\date{}
\maketitle
\raggedbottom
%\date{}

%%%%%%%%%%%%%%%%%%%%%%%%%%%%%%

%%%%%%%%%%%%%%%%%%%%%%%%%%%%%%%%%%%%%%%%%%%%%%%%%%%%%%%%
%%%%%%%%%%%%%%%%%%%%%%%%%%%%%%%%
%\begin{center}
%\noindent
%{\bf  Abstract}~
%\end{center}

\begin{abstract}
In rough stochastic PDE theory of Hairer type,
rough path lifts with respect to the space variable
of two-parameter continuous Gaussian processes 
play a main role.
A prominent example of such processes is the solution of the stochastic heat equation 
under the periodic condition.
The main objective of this paper is to show that the law of the spatial lift of this process 
satisfies a Schilder type large  deviation principle
on the continuous path space over a geometric rough path space.
\end{abstract}

%\noindent
%{\bf  Key words and phrases}~
%%%%%%%%%%%%%%%%%%%%%%%%%%%

%%%%%%%%%%%%%%%%%%%%%%%%%%%%%%%

\section{Introduction and main result}
In rough path theory of T. Lyons, 
the notion of paths is generalized to a great extent
and so is that of ordinary differential equations.
They are called rough paths and rough differential equations (RDEs), respectively.
The solution map of an RDE is called an It\^o map,
which is defined for every rough path and, moreover, is continuous 
with respect to the topology of rough path space (Lyons' continuity theorem).
As a result, stochastic differential equations (SDEs) in the usual sense 
are made deterministic or "dis-randomized".

Even though It\^o maps are deterministic, 
the probabilistic aspect of the theory is still very important undoubtedly.
In a biased view of the author, 
a large deviation principle of Schilder type is a central issue
in stochastic analysis on rough path spaces.
This kind of large deviations was first shown by 
Ledoux, Qian, and Zhang \cite{lqz} for the law of  Brownian rough path.
Combined with Lyons' continuity theorem, 
this result immediately recovers well-known Freidlin-Wentzell type 
large deviations for solutions of SDEs.
Since then many papers have been published on this topic \cite{der,dd, fv05, fv07, ina, ik, ms}.

Naturally, one would like to apply rough path theory
to stochastic PDEs.
There have been some successful attempts.
In this paper, we focus on M. Hairer's theory \cite{hai1, hw, hmw},
which is based on M. Gubinelli's controlled path formalism.
In this theory, rough path integration is used for the space variable $x \in S^1 ={\bf R}/{\bf Z}$
for each fixed time variable $t>0$.
This is surprising because almost everyone regarded 
solutions of stochastic PDEs as  processes indexed by $t$-variable that take values in 
function spaces of $x$-variable
and then  modify and apply infinite dimensional rough path theory.
Not only this point of view is novel, 
but this theory also turned out to be very powerful
when he rigorously solved KPZ equation in the periodic case for the first time \cite{hai2}.

Under these circumstances, it seems natural and necessary 
to develop stochastic analysis in this framework.
In this paper we will prove a large deviation principle 
of Schilder type for the spatial lift of the (scaled) solution $\psi$ of 
the stochastic heat equation on $S^1$.
This process $\psi$ plays a crucial role in \cite{hw, hmw}.
To our knowledge, 
a large deviation principle is new in rough stochastic PDE theories of any kind.

Now we introduce our setting. 
We will give precise definitions and 
detailed explanations in later sections.
Let us recall the stochastic heat equation on $S^1$.
As usual $S^1 ={\bf R}/{\bf Z}$ is regarded as $[0,1]$
with the two end points identified
and $\triangle = \triangle_{S^1}$ stands for the periodic Laplacian.
Let 
$\xi^i =\xi (t,x)^i ~(1 \le i \le d)$  
are independent copies of  the space-time white noise 
associated with $L^2([0,T] \times S^1)$
with the (formal) covariance 
${\mathbb E}[ \xi (t,x)^i \xi (s,y)^j]= \delta_{ij}\cdot\delta_{t-s}\cdot \delta_{x-y}$.
Let $\psi= \psi(t,x)$ be a unique solution of the following ${\bf R}^d$-valued stochastic PDE.
\begin{equation}\label{eq.psi.spde}
\partial_t \psi = \triangle_x \psi + \xi, \qquad \mbox{with} \qquad \psi(0, x) \equiv 0.
\end{equation}
Then,
$\psi =(\psi(t,x))_{0 \le t \le T, 0 \le x \le 1}$ 
is a two-parameter continuous Gaussian process.
It was shown in \cite{hw} that, 
(i) for each $t$, $x \mapsto \psi(t,x)$ admits a natural lift 
to a geometric rough path $(x,y) \mapsto \Psi (t; x,y)$ a.s.
and (ii) there exists a  modification of $\Psi$ such that
$t \mapsto \Psi (t; \bullet,\star)$ is continuous in the geometric rough path space a.s. 
In this theory, 
a solution of a rough stochastic PDE is obtained from $\Psi$.
Therefore, it is important to analyze (the law of) $\Psi$.

Let $1/3 <\a <1/2$.
We denote by $G \Omega^H_{\a} ({\bf R}^d)$
the $\a$-H\"older geometric rough path space over ${\bf R}^d$.
The first level path of $X \in G \Omega^H_{\a} ({\bf R}^d)$
is a usual path in ${\bf R}^d$ which starts at $0$.
Let $G \hat\Omega^H_{\a} ({\bf R}^d) \cong {\bf R}^d \times G \Omega^H_{\a} ({\bf R}^d)$
be the $\a$-H\"older geometric rough path space
in an extended sense so that information of 
 the initial values of the first level paths are added.
For each $t$, the random variable
$\Psi (t; \bullet,\star)$ takes values in this Polish space
$G \hat\Omega^H_{\a} ({\bf R}^d)$.
Let ${\cal P}_{\infty} G \hat\Omega^H_{\a} ({\bf R}^d) 
= C([0,T], G \hat\Omega^H_{\a} ({\bf R}^d))$
be the continuous path space over $G \hat\Omega^H_{\a} ({\bf R}^d)$.
Its topology is given by the uniform convergence in $t$ as usual.
The random variable $\Psi$ takes values in this Polish space
and hence its law  is a probability measure on this space.

Introduce a small parameter $0 < \ve \le 1$.
Let $\ve \Psi$ is the dilation of $\Psi$ by $\ve$, which is equal to
the natural lift of $\ve \psi$, anyway.
Denote by $\nu_{\ve}$
 the law of $\ve \Psi$ on ${\cal P}_{\infty} G \hat\Omega^H_{\a} ({\bf R}^d)$.
Our main result is the following:

\vspace{3mm}

\noindent
{\bf Main result:}~
{\it For any $\a \in (1/3, 1/2)$, the family $(\nu_{\ve})_{0 <\ve \le 1}$ 
of probability measures on ${\cal P}_{\infty} G \hat\Omega^H_{\a} ({\bf R}^d)$
satisfies a large deviation principle as $\ve \searrow 0$ with a good rate function $I$.}

\vspace{3mm}

\noindent
See Theorem \ref{pr.ldp.st} below for a precise statement of our main result.
A concrete expression of $I$ will be given in Section \ref{sec.ldp}.
A comment on a Freidlin-Wentzell type large deviation principle 
for solutions of rough stochastic PDEs will also be given after the main theorem.

The organization of this paper is as follows.
In Section \ref{sec.grps}
we introduce  several kind of path spaces over geometric rough path spaces
and show basic properties between them.
In Section \ref{sec.cov} we prove a few properties of 
 the covariance of $\psi$ in one-dimensional case.
We use them in the following section to prove convergence of the dyadic polygonal approximations.

In Section \ref{sec.dyconv}
we prove that the spatial 
lift of the dyadic polygonal approximations of $\psi$
converge with respect to a Besov type topology which is stronger than the topology of
${\cal P}_{\infty} G \hat\Omega^H_{\a} ({\bf R}^d)$.
Unlike preceding examples in \cite{fggr, hai1, hmw, hw, r},
 in which Kolmogorov's continuity criterion is used to obtain regularity in $t$-variable, 
we directly prove a.s. and $L^p$-convergence in function spaces of space-time variables.
In this respect, the argument in this section may be new.

In Section \ref{sec.ldp} we prove our main result, 
namely, large deviations for the law of $\ve \Psi$.
Our method is based on Friz and Victoir's in \cite{fv07}.
%Besides its generality, 
%
Their method is general and moreover is available
%an advantage of their method is its availability 
even when  regularity of Cameron-Martin paths is not understood very well.
 (Sometimes, lifting a Gaussian process 
 may be easier than lifting its Cameron-Martin paths.)
Neither
in our case do we know much about (the lift of) Cameron-Martin elements of $\psi$. 
So we use their method.

Throughout this paper, $c$ denotes an unimportant positive constant  
which may change from line to line.

%\newpage
%%%%%%%%%%%%%%%%%%%%%%%%%%%%%%%%%%%%%%%%%%%%%%%%%%%%%%%%%%%%%%%%%%%%%%
%%%%%%%%%%%%%%%%%%%%%%%%%%%%%%%%%%%%%%%%%%%%%%%%%%%%%%%%%%%%%%%%%%%%%
%%    New  Section 
%%%%%%%%%%%%%%%%%%%%%%%%%%%%%%%%%%%%%%%%%%%%%%%%%%%%%%%%%%%%%%%%%%%%%%
%%%%%%%%%%%%%%%%%%%%%%%%%%%%%%%%%%%%%%%%%%%%%%%%%%%%%%%%%%%%%%%%%%%%%%

\section{Path space over geometric rough path space}\label{sec.grps}

In this section we recall the definitions and some basic facts of 
the geometric rough path spaces, which are endowed with
H\"older, Besov,  and the uniform norms.
Then we consider continuous path spaces over geometric rough path spaces
and introduce several kinds of norms on them.
All the ingredients of this section are either known or easily derived from known facts.
We refer to Appendix A.2 in Friz and Victoir \cite{fvbk}
for basic information on H\"older and Besov norms.

%%%%%%

%%%%%%%%%%%%%%%%%%%%%%%%%

%\vspace{10mm}
%\noindent
%$\spadesuit$~
%
%
Let us first recall 
H\"older and Besov norms on continuous path space over a real Banach space.
Let ${\cal V}$ be a real Banach space. 
The set of continuous path $f:[0,T] \to {\cal V}$ is denoted by $C([0,T], {\cal V})$.
It is a real Banach space with the usual uniform norm $\|f\|_{\infty} :=\sup_{0 \le t \le T} |f_t|_{{\cal V}}$.
The subspace of all the continuous paths which start at $0$ is denoted by $C_0 ([0,T], {\cal V})$.

For $\a \in (0,1]$ and $f \in C([0,T], {\cal V})$, we define $\a$-H\"older norm of $f$ by
\[
\|f\|_{H ;\a} := |f_0 |_{{\cal V}} 
+
\sup_{0 \le s< t \le T} \frac {|f_t -f_s|_{{\cal V}} }{|t-s|^{\a}}.
\]
The subspace of all the paths with $\|f\|_{H ;\a} <\infty$
is denoted by $C^{H;\a} ([0,T], {\cal V})$.
Next we introduce Besov norm.
For $\a \in (0,1]$ and $m \ge 1$ with $\a >1/m$, we define $(\a,m)$-Besov norm of $f$ by
\[
\|f\|_{B ;\a,m} := |f_0 |_{{\cal V}} 
+
\Bigl\{
\iint_{ {\cal S}(T)}  \frac {|f_t -f_s|^{m}_{{\cal V}} }{|t-s|^{1 +\a m}} dsdt
\Bigr\}^{1/m}.
\]
Here, ${\cal S}(T)= \{(s,t)  ~|~0 \le s \le t \le T \}$.
(When $T=1$ we write ${\cal S}={\cal S}(1)$ for simplicity).
The subspace of all the paths with $\|f\|_{B ;\a,m} <\infty$
is denoted by $C^{B;\a,m} ([0,T], {\cal V})$.
It is obvious that $\|f\|_{B ;\a,m} \le c \|f\|_{H ;\a'}$
for some constant $c=c_{\a,\a',m} >0$ if $\a <\a'$.
By Besov-H\"older embedding theorem (e.g. Corollary A.2, \cite{fvbk}), 
it also holds that 
$\|f\|_{H ;\a - 1/m} \le c' \|f\|_{B ;\a,m}$
for some constant $c'=c'_{\a,m} >0$.
In particular, $C^{B;\a,m} ([0,T], {\cal V})$ is continuously imbedded in $C^{H;\a -1/m} ([0,T], {\cal V})$.
The subspace of all the Besov or H\"older paths which start at $0$
is denoted with the subscript "$0$".

%\vspace{10mm}\noindent
%$\spadesuit$~ 

Let us next recall 
geometric rough path spaces.
%
%
%\noindent
Throughout this paper, the parameter interval for rough paths
is $[0,1]$ and we set 
${\cal S} ={\cal S} (1)=\{ (x,y) ~|~ 0 \le x \le y \le 1 \}$.
Let $A \in C( {\cal S},  {\bf R}^d)$ vanish on the diagonal.
(The totality of such $A$'s will be denoted by $C_0( {\cal S},  {\bf R}^d)$.)
We set 
\begin{align}
\|A\|_{ H;\alpha }  &= \sup_{0 \le x < y \le 1  } \frac{| A_{x,y} | }{ |y-x|^{\alpha}}
\qquad
&
(0< \alpha \le 1),
\label{hld.def}
\\
\| A\|_{B;\alpha, m}  &= 
\Bigl(
\iint_{ {\cal S} } 
 \frac{ | A_{x,y} |^m }{  |y-x|^{1 +m\alpha}} dxdy  \Bigr)^{1/m}
\qquad
&
( 0< 1/m <\alpha \le 1).
\label{bsv.def}
\end{align}
These are called $\alpha$-H\"older norm and $(\alpha, m)$-Besov norm, respectively.

Let $T^2 ({\bf R}^d) = {\bf R} \oplus {\bf R}^d \oplus ({\bf R}^d)^{\otimes 2}$
be the truncated tensor algebra of step 2.
The set of elements of the form $(1, \bullet, \star)$ forms a non-abelian group 
under the tensor multiplication $\otimes$.
The unit element is ${\bf 1}=(1, 0,0)$.
A scalar action $(1, a_1, a_2) \mapsto (1, \lambda a_1, \lambda^2 a_2)$ for $\lambda \in {\bf R}$
is called the dilation.

A continuous map
$A=(1, A^1, A^2) :{\cal S}\to   T^2 ({\bf R}^d)$
is called multiplicative if it satisfies that
\begin{align}
A^1_{x,y} = A^1_{x,z}+A^1_{z,y}, 
 \qquad
  A^2_{x,y} = A^2_{x,z}+A^2_{z,y} +  A^1_{x,z} \otimes A^1_{z,y},
\qquad
(x \le z \le y).
\label{chen.eq}
\end{align}
This relation is called Chen's identity and 
can also be written as $A_{x,y} = A_{x,z} \otimes A_{z,y}$, 
where $\otimes$ stands for the multiplication of $T^2 ({\bf R}^d)$.
In particular, $A_{x,y}$ is a "difference" of a group-valued path, since 
$A_{x,y} = (A_{0,x})^{-1} \otimes A_{0,y}$.
Note that $A^1$ and $A^2$ vanish on the diagonal if they are multiplicative. 
We denote by $\Omega_{\infty} ({\bf R}^d)$ the set of such continuous multiplicative functionals.
A distance on $\Omega_{\infty} ({\bf R}^d)$ is given by 
$d(A, B) = \|A^1 - B^1\|_{ \infty} + \|A^2 -B^2 \|_{\infty}$,
where $\|\,\cdot\,\|_{ \infty}$ denotes the sup-norm over ${\cal S}$ as usual.

Let $1/3 <\alpha \le 1/2$.
The space of ${\bf R}^d$-valued $\alpha$-H\"older rough path is defined by
\begin{align}
\Omega_{\alpha}^H ({\bf R}^d) 
&= \{  A=(1, A^1, A^2) \in C( {\cal S}, T^2( {\bf R}^d) )
\nn
\\ 
&  \qquad 
~|~  \mbox{multiplicative and   }  \|A^1 \|_{ H;\alpha } <\infty,  \|A^2 \|_{H;2 \alpha} <\infty\}.
\nn
\end{align}
The topology of this space is naturally induced by the following distance:
$d(A, B) = \|A^1 - B^1\|_{ H;\alpha } + \|A^2 -B^2 \|_{H;2 \alpha }$.
In the same way, $(\alpha, m)$-Besov  rough path is 
defined for $m \ge 2$ and $1/3 <\alpha \le 1/2$ with $\alpha -m^{-1}> 1/3$ as follows;
\begin{align}
\Omega_{m, \alpha}^B   ({\bf R}^d) 
&= \{  A=(1, A^1, A^2) \in C( {\cal S}, T^2( {\bf R}^d) )
\nn
\\ 
&  \qquad 
~|~  \mbox{multiplicative and   }  \|A^1 \|_{B;\alpha, m} <\infty, 
 \|A^2 \|_{B; 2\alpha, m/2 } <\infty\}.
\nn
\end{align}
The topology of this space is naturally induced by the following distance:
$d(A, B) = \|A^1 - B^1\|_{B;\alpha, m } + \|A^2 -B^2 \|_{ B; 2 \alpha, m/2}$.
In what follows, we will often write $A =(A^1, A^2)$ for simplicity, 
since the $0$th component "$1$" is obvious.

A Lipschitz continuous path (i.e., $1$-H\"older continuous path) 
$x \in C_0^{H;1} ( [0,1], {\bf R}^d)$
admits a natural lift to a rough path by setting 
\[
A^1_{x,y} := a_y -a_x, \qquad A^2_{x,y} := \int_x^y (a_z -a_x) \otimes da_z,  \qquad (x,y) \in {\cal S}. 
\]
It is easy to see that $A \in \Omega_{\alpha}^H ({\bf R}^d) \cap \Omega_{\alpha, m}^B ({\bf R}^d)$.
We call a rough path $A$ obtained in this way a smooth rough path lying above $a$,
or the natural lift of $a$.
The natural lift map is denoted by ${\cal L}_2$, i.e., $A ={\cal L}_2 (a)$.

Now we introduce geometric rough path spaces.
Let $G\Omega_{\alpha}^H ({\bf R}^d)$ be the closure 
of the set of smooth rough paths in $\Omega_{\alpha}^H ({\bf R}^d)$. 
This is called the geometric rough path space with $\alpha$-H\"older norm.
The geometric rough path space $G \Omega_{ \alpha, m}^B   ({\bf R}^d) $ 
with $(\alpha, m)$-Besov norm is similarly defined.
Hence, we have the following inclusions;
\[
G\Omega_{\alpha}^H  ({\bf R}^d)  \subset 
\Omega_{ \alpha}^H   ({\bf R}^d),
\qquad
G\Omega_{\alpha, m}^B   ({\bf R}^d)  \subset 
\Omega_{\alpha, m}^B   ({\bf R}^d).
\]

Moreover, by Besov-H\"older embedding in Proposition \ref{inj.HB.pr} below, 
we also have 
\[
G\Omega_{\alpha,m}^B  ({\bf R}^d)  \subset 
G\Omega_{ \alpha -1/m}^H   ({\bf R}^d),
\qquad
\Omega_{\alpha, m}^B   ({\bf R}^d)  \subset 
\Omega_{\alpha -1/m}^H   ({\bf R}^d).
\]
if $1/3< \a \le 1/2$, $m \ge 2$, and $\a -1/m >1/3$.
Note that these continuous inclusions are bounded 
(in the sense that they map any bounded set to a bounded set).

\begin{prop}\label{inj.HB.pr}
Assume $1/3< \a \le 1/2$, $m \ge 2$, and $\a -1/m >1/3$.
Then, $\Omega_{\alpha, m}^B  ({\bf R}^d)$ is continuously embedded 
in $\Omega_{\alpha -1/m}^H  ({\bf R}^d)$.
Consequently, 
{\rm (i)}~$G\Omega_{\alpha, m}^B  ({\bf R}^d)$ is continuously embedded 
in $G\Omega_{\alpha -1/m}^H  ({\bf R}^d)$
and 
{\rm (ii)}~$\Omega_{\alpha, m}^B  ({\bf R}^d)$ is a complete metric space
and 
$G\Omega_{\alpha, m}^B  ({\bf R}^d)$ is a Polish space.
\end{prop}

\Proof
Recall the following inequalities for Besov-H\"older embedding: 
For some positive constant $c=c_{\a,m}$ independent of $A, B \in \Omega_{\infty}  ({\bf R}^d)$, 
we have
\begin{align}
\| A^1 -B^1\|_{H; \a -1/m } &\le c \| A^1 -B^1\|_{B;\a,m }, 
\nn\\
 \| A^2 -B^2\|_{H; 2\a -2/m } &\le 
 c \bigl( \| A^1 -B^1\|_{B;\a,m }+\| A^2 -B^2\|_{B;2\a,m/2 }\bigr)
 \nn\\
  & \quad \times
   \bigl( \| A^1\|_{B;\a,m }+\| A^2\|_{B;2\a,m/2 } + \| B^1\|_{B;\a,m }+\| B^2\|_{B;2\a,m/2 }\bigr).
    \label{est_FV.A10}
\end{align}
The first one has already been explained.
The second one is found in Proposition A.10, pp. 576--579, \cite{fvbk}.
(In the proof, only multiplicativity of $A$ and $B$ is used. 
In other words, $t \mapsto A_{0,t}, B_{0,t}$ need not take their values in 
the free nilpotent group of step 2. See \cite{fvbk} for details.)
From these inequalities, we can easily see that 
$\Omega_{\alpha, m}^B  ({\bf R}^d)$ is continuously embedded 
in $\Omega_{\alpha -1/m}^H  ({\bf R}^d)$.
%%%

Now we prove the rest of the proposition.
The only non-trivial part is  completeness of $\Omega_{\alpha, m}^B  ({\bf R}^d)$.
Let $A(n)~(n=1,2,\ldots)$ be a Cauchy sequence in $(\a, m)$-Besov topology.
In the Besov topology, there exists a limit $A(\infty)$. 
But, 
$A(\infty)$ is continuous and multiplicative on ${\cal S}$
since convergence is also in $(\a -1/m)$-H\"older topology.
Hence, $\Omega_{\alpha, m}^B  ({\bf R}^d)$ is complete.
\QED

%%%%%%%%%%%%%%%%%%%%%%%%%%%%%%%%%%%%%%%%%%%%%%%%%
%\vspace{10mm}
%\noindent
%$\spadesuit$~Initial value.

In the above definition, the first level path of a rough path 
is naturally identified with a path in the usual sense which starts at origin.
Now we slightly modify the definition so that 
the first level path can start at any point.

Set $\hat\Omega_{\infty}  ({\bf R}^d) = {\bf R}^d \times \Omega_{\infty}  ({\bf R}^d)$.
The distance on it is the natural one for a product space.
The path $[0,1 ] \ni x \mapsto v + A^1_{0,x}$ is said to be 
the first level path of  $(v, A) \in \hat\Omega_{\infty}  ({\bf R}^d)$.
$v$ is called the initial value
The dilation 
naturally extends on this space by $\lambda  (v, A) = (\lambda v, \lambda  A)$.
In the same way, $\hat\Omega_{\alpha, m}^B  ({\bf R}^d)$ and
 $\hat\Omega_{\alpha}^H  ({\bf R}^d)$ are defined.

Conversely, for $a \in C^{H;1}([0,1], {\bf R}^d)$,
we set ${\cal L}_2 (a) = (a_0, {\cal L}_2(a_{\cdot} -a_0) ) \in \hat\Omega_{\infty}  ({\bf R}^d)$
and call it a natural lift of $a$ or a smooth rough path lying above $a$.
Geometric rough path spaces $G\hat\Omega_{\infty}  ({\bf R}^d)$, 
$G\hat\Omega_{\alpha, m}^B  ({\bf R}^d)$, 
and $G\hat\Omega_{\alpha}^H  ({\bf R}^d)$ in this extended sense 
are defined as the closure of the set of smooth rough paths 
as before.

\begin{remark} \label{rem.i.val}
{\rm (i)}~For $(v, A) =(v, A^1,A^2) \in \hat\Omega_{\infty}  ({\bf R}^d)$,
we will sometimes write $(a,A^2)$, where $a_x = v+ A^1_{0,x}$.
In particular, for $b \in C^{H;1}([0,1], {\bf R}^d)$, 
$(b,B^2)$ stands for
the smooth rough path lying above $b$.
\\
\noindent
{\rm (ii)}~ Almost all the results for $G\Omega_{\alpha}^H  ({\bf R}^d)$ etc. 
also hold for $G\hat\Omega_{\alpha}^H  ({\bf R}^d)$ etc. 
with trivial modifications.
\\
\noindent
{\rm (iii)}~$\Omega_{\alpha}^H  ({\bf R}^d)$ is basically the same 
rough path space as in \cite{hai1, hw, hmw, hai2}, etc.,
in which the rough path space is denoted by ${\cal D}^{\a} ({\bf R}^d)$.
The only difference is that 
in this paper a rough path is defined on the simplex ${\cal S}$,
while in these papers it is defined on the simplex $[0,1]^2$.
However, under Chen's identity  values of rough paths
on $[0,1]^2 \setminus {\cal S}$ are automatically determined.
Hence, there is essentially no difference.
\end{remark}

%%%%%%%%%%%%%%%%%%%%%%%%%%%%%%%%%%%%%%%%%%%%%%%%%
%\vspace{10mm}
%\noindent
%$\spadesuit$~

Now we consider several kind of path spaces over geometric rough path space.
Let $T>0$ and assume $1/3 <\a \le 1/2$, $m \ge 2$, and $\a -1/m >1/3$.
We first define continuous path spaces over geometric rough path spaces
with the usual $\sup$-distance.
Let ${\cal P}_{\infty} G\hat\Omega_{\alpha}^H  ({\bf R}^d) = C([0,T], G\hat\Omega_{\alpha}^H  ({\bf R}^d))$ 
be a continuous path space 
over $G\hat\Omega_{\alpha}^H  ({\bf R}^d)$.
An element of this set is of the form $(v_t, A^1(t;x,y), A^2(t;x,y))$.
(We often write $A^i(t;x,y) = A^i_t(x,y)$ for simplicity.)
The distance on this space is defined by
\begin{align}
{\rm dist} \bigl( (v, A^1, A^2), (u, B^1, B^2)\bigr)
&=
\sup_{t \in [0,T]} |v_t -u_t|
\nn\\
&+
\sup_{t \in [0,T]} \|A^1_t -B^1_t \|_{H;\a}
+
\sup_{t \in [0,T]} \|A^2_t -B^2_t \|_{H; 2\a}.
\nn
\end{align}
This space is the one used in M. Hairer's rough stochastic PDE theory (see \cite{hai1, hw, hmw, hai2}).
In the same way, we define
${\cal P}_{\infty} G\hat\Omega_{\alpha,m}^B  ({\bf R}^d)$
and 
${\cal P}_{\infty} G\hat\Omega_{\infty}  ({\bf R}^d)$
as well as their distances.

Let $0<\b \le 1$.
Set 
\begin{align}
{\cal P}_{\b}^H G\hat\Omega_{\alpha}^H  ({\bf R}^d)
&=
\{ (v, A^1, A^2) \in {\cal P}_{\infty} G\hat\Omega_{\alpha}^H  ({\bf R}^d)~|~
\nn\\
& \qquad 
\mbox{$t \mapsto v_t, A^1_t$ is $\b$-H\"older continuous,
$A^2_t$ is $2\b$-H\"older continuous}
\}.
\nn
\end{align}
Of course, $t \mapsto A^i_t$ is $i\b$-H\"older continuous with respect to $\| \,\cdot\,\|_{H; i\a}$
for $i=1,2$,
that is, 
$\sup_{0 \le s <t \le T}  \| A^i_t -A^i_s \|_{H; i\a}/(t-s)^{i\b}<\infty$.
The distance on this space is defined by
\begin{align}
\lefteqn{
{\rm dist}  \bigl( (v, A^1, A^2), (u, B^1, B^2)\bigr)
}
\nn\\
&=
 \|v -u\|_{H;\b}
+
\sum_{i=1,2}
{\cal N} \bigl[
A^i -B^i ; C^{H;i\b}([0,T], C_0^{H;i\a}({\cal S},  ({\bf R}^d)^{ \otimes i}))
\bigr]
\nn\\
&= \|v -u\|_{H;\b}
+\sum_{i=1,2} \| A^i_0- B_0^i \|_{H;i\a}
\nn\\
& \qquad +
\sum_{i=1,2}
\sup_{s<t, x<y} 
\frac{ | (A^i_t(x,y) - B^i_t(x,y)) -(A^i_s(x,y) - B^i_s(x,y) ) | }{ (t-s)^{i\b} (y-x)^{i\a} }.
\nn
\end{align}
Here ${\cal N}[\,\cdot\,;{\cal V}]$ denotes the norm of a Banach space ${\cal V}$.
In the same way, ${\cal P}_{\b}^H G\hat\Omega_{\alpha,m}^B  ({\bf R}^d)$ and 
${\cal P}_{\b}^H G\hat\Omega_{\infty}  ({\bf R}^d)$ can be defined,
but details are omitted.

Assume in addition that $\b >1/m$.
We set $(\b,m)$-Besov path spaces over geometric rough path spaces as follows;
\begin{align}
&{\cal P}_{\b,m}^B G\hat\Omega_{\alpha,m}^B  ({\bf R}^d)
=
\bigl\{ (v, A^1, A^2) \in {\cal P}_{\infty} G\hat\Omega_{\alpha,m}^B  ({\bf R}^d)~|~
\nn\\
&\mbox{$t \mapsto v_t, A^1_t$ is $(\b,m)$-Besov continuous,
$A^2_t$ is $(2\b,m/2)$-Besov continuous} \bigr\}.
\nn
\end{align}
Of course, $t \mapsto A^i_t$ is $(i\b, m/i)$-Besov continuous with respect to $\| \,\cdot\,\|_{B; i\a, m'/i}$
for $i=1,2$.
The spaces
${\cal P}_{\b,m}^B G\hat\Omega_{\alpha}^H ({\bf R}^d)$ and 
${\cal P}_{\b,m}^B G\hat\Omega_{\infty} ({\bf R}^d)$,
and the distances on these spaces are also defined in a similar way.
The distance on ${\cal P}_{\b,m}^B G\hat\Omega_{\alpha,m}^B  ({\bf R}^d)$
has the following form;
\begin{align}
\lefteqn{
{\rm dist} \bigl( (v, A^1, A^2), (u, B^1, B^2)\bigr)
}
\nn\\
&=
 \|v -u\|_{B;\b,m}
+
\sum_{i=1,2}
{\cal N} \bigl[
A^i -B^i ; C^{B;i\b,m/i}([0,T], C_0^{B;i\a, m/i}({\cal S},  ({\bf R}^d)^{ \otimes i}))
\bigr]
\nn\\
&= \|v -u\|_{B;\b,m}
+\sum_{i=1,2} \| A^i_0- B_0^i \|_{B;i\a, m/i}
\nn\\
& \qquad +
\sum_{i=1,2}
\Bigl\{
\iint_{{\cal S}(T)}dsdt
\iint_{{\cal S}}
\frac{| (A^i_t(x,y) - B^i_t(x,y)) -(A^i_s(x,y) - B^i_s(x,y) ) |^{m/i }}
{(t-s)^{1+\b m} (y-x)^{1+ \a m}}
dxdy
\Bigr\}^{i/m}.
\nn
\end{align}
Finally, note that all of these ${\cal P} G\hat\Omega$'s introduced above are complete.

%%%%%%%%%%%%%%%%%%%%%%%%%%%%%%%%%%%%%%%%%%%%%%%%%
%\vspace{10mm}
%\noindent
%$\spadesuit$~Inclusion of path spaces.
%
There are of course natural inclusions between these ${\cal P} G\hat\Omega$'s.
Now we discuss two of them for later use.
\begin{prop}\label{pr.PGL.inj}
Let $1/3< \a' <\a \le 1/2$, $0 <\b \le 1$, $m \ge 2$ such that 
$\a' < \a -1/m$ and $\b >1/m$.
Then, we have the following bounded, continuous inclusions;
\[
{\cal P}_{\b,m}^B G\hat\Omega_{\alpha,m}^B  ({\bf R}^d)
\hookrightarrow
{\cal P}_{\infty} G\hat\Omega_{\alpha'}^H  ({\bf R}^d) 
\hookrightarrow
{\cal P}_{\infty} G\hat\Omega_{\infty}  ({\bf R}^d). 
\]
Moreover, the left inclusion is compact in the sense that it maps any bounded subset 
to a precompact subset.
\end{prop}

\Proof
The right inclusion is obvious.
Now we consider the left one.
From (\ref{est_FV.A10}) and Besov-H\"older embedding in $t$-variable,
 we easily obtain
\begin{align}
\sup_t \| A^1_t -B^1_t\|_{H; \a'} &\le c \| A^1 -B^1\|_{B;\b, \a,m }, 
\nn\\
 \sup_t  \| A^2_t -B^2_t\|_{H; 2\a' } &\le 
 c \bigl( \| A^1 -B^1\|_{B;\b,\a,m }+\| A^2 -B^2\|_{B;2\b,2\a,m/2 }\bigr)
 \nn\\
   \quad \times
   \bigl( \| A^1\|_{B;\b, \a,m } &+\| A^2\|_{B;2\b, 2\a,m/2 }
       + \| B^1\|_{B;\b, \a,m }+\| B^2\|_{B;2\b, 2\a,m/2 }\bigr).
    \label{est_FV.A++}
\end{align}
Here, the norms stand for 
${\cal N} \bigl[ \,\cdot\,; C^{B;i\b,m/2}([0,T], C_0^{B;i\a, m/i}({\cal S},  ({\bf R}^d)^{ \otimes i}))
\bigr]$ for $i=1,2$.
Inequalities (\ref{est_FV.A++}) show that the left map is a continuous inclusion.

We prove compactness.
First, the inclusion $G\hat\Omega_{\alpha,m}^B  ({\bf R}^d)
\hookrightarrow
G\hat\Omega_{\alpha'}^H  ({\bf R}^d)$ is compact.
Second, if $\{ (v(n), A(n)^1, A(n)^2) \}_{n=1,2,\ldots}$ is a bounded sequence 
in ${\cal P}_{\b,m}^B G\hat\Omega_{\alpha,m}^B  ({\bf R}^d)$,  
they are bounded and uniformly continuous $G\hat\Omega_{\alpha,m}^B  ({\bf R}^d)$-valued paths
 since $\beta >1/m$ and Besov-H\"older embedding.
By (\ref{est_FV.A10}) they are also bounded and uniformly continuous
as $G\hat\Omega_{\alpha'}^H  ({\bf R}^d)$-valued paths.

Now we can use Ascoli-Arzela type argument as follows.
By the diagonalization argument, 
there exists a subsequence $\{n_k\}_{k=1,2,\ldots}$
such that
$\{ (v(n_k)_t, A(n_k)^1_t, A(n_k)^2_t) \}_{k=1,2,\ldots}$
converges in $G\hat\Omega_{\alpha'}^H  ({\bf R}^d)$
for any $t \in [0,T] \cap {\mathbf Q}$.
Let $\ve >0$ be arbitrary.
By the uniform continuity and compactness of $[0,T]$, 
there exists finitely many (relatively) open intervals ${\cal U}_j~(1 \le j \le l)$
such that {\rm (i)}~ $[0,T] = \cup_{j=1}^l {\cal U}_j$
and {\rm (ii)}~it satisfy that 
$$
{\rm d} \bigl( (v(n)_t, A(n)^1_t, A(n)^2_t),(v(n)_s, A(n)^1_s, A(n)^2_s)\bigr)  <\ve
\qquad
(s,t \in {\cal U}_j, \, 1 \le j \le l).
$$
Here, ${\rm d}$ denotes the distance on $G\hat\Omega_{\alpha'}^H  ({\bf R}^d)$.
Next
choose $t_j \in {\cal U}_j \cap {\bf Q} \, (1 \le j \le l)$.
Since $l$ is finite, there exists $k_0 \in {\bf N}$ such that 
if $k, m \ge k_0$ then
$$
{\rm d} \bigl( (v(n_k)_{t_j}, A(n_k)^1_{t_j}, A(n_k)^2_{t_j}),
(v(n_m)_{t_j}, A(n_m)^1_{t_j}, A(n_m)^2_{t_j})\bigr)  <\ve
\qquad
( 1 \le j \le l).
$$
The $3\ve$-argument implies that if $k, m \ge k_0$ then
$$
\sup_{0 \le t \le T} {\rm d} \bigl( (v(n_k)_{t}, A(n_k)^1_{t}, A(n_k)^2_{t}),
(v(n_m)_{t}, A(n_m)^1_{t}, A(n_m)^2_{t})\bigr)  < 3\ve,
$$
which shows that this subsequence is Cauchy 
in ${\cal P}_{\infty} G\hat\Omega_{\alpha'}^H  ({\bf R}^d)$.
\QED

%%%%%%%%%%%%%%%%%%%%%%%%%%%%%%%%%%%%%%%%%%%%%%%%%%%%%%%%%%%%%%%%%%%%%%%%%%%%%%%%%%%%%
%%%%%%%%%%%%%%%%%%%%%%%%%%%%%%%%%%%%%%%%%%%%%%%%%%%%%%%%%%%%%%%%%%%%%%%%%%%%%%%%%%%%%
%   New  SECTION
%%%%%%%%%%%%%%%%%%%%%%%%%%%%%%%%%%%%%%%%%%%%%%%%%%%%%%%%%%%%%%%%%%%%%%%%%%%%%%%%%%%%
%%%%%%%%%%%%%%%%%%%%%%%%%%%%%%%%%%%%%%%%%%%%%%%%%%%%%%%%%%%%%%%%%%%%%%%%%%%%%%%%%%%%%%%%%\newpage

\section{Covariance }\label{sec.cov}
The most important data for Gaussian processes are their covariances.
In this section we will calculate the covariance of $\psi$.
Throughout this section, the dimension is always $d=1$.

Let $t \ge 0$ and $x \in {\mathbf R}/{\mathbf Z} \cong S^1$.
As usual, we identify $S^1$ with $[0,1]$.
$\triangle = \triangle_{S^1}$ denotes the periodic Laplacian. 
$\xi =\xi (t,x)$ is the space-time white noise 
with the formal covariance ${\mathbb E}[ \xi (t,x) \xi (s,y)]= \delta_{t-s}\cdot \delta_{x-y}$.
Let $\psi= \psi(t,x)$ be a unique solution of the following real-valued stochastic PDE.
\[
\partial_t \psi = \triangle_x \psi + \xi, \qquad \mbox{with} \qquad \psi(0, x) \equiv 0.
\]
$\psi$ is a two-parameter Gaussian process and can be written down as follows;
\begin{align}
\psi (t,x)  
=
\Bigl[ \int_0^t  e^{(t-s) \triangle} \xi(s, \cdot) ds  \Bigr] (x)
%\nn\\
=
 \int_0^t ds \int^1_0 \tilde{p}_{t-s} (x-y) \xi (s,y) dy.
\nn
 \end{align}
Here, 
\[
\tilde{p}_t (x)=  \sum_{n = -\infty}^{\infty}  p_t (x+n),
\qquad
\quad
p_t (x)=\frac{1}{\sqrt{4 \pi t}} \exp \bigl(-\frac{x^2 }{4t} \Bigr)
\]
are the heat kernels of $S^1$ and ${\mathbf R}$, respectively.

%%%%%
First, we give an explicit expression of the covariance.
Observe that it has a few kinds of symmetries.

%%%%%%%%%%%%%%%%%%%%%%
\begin{lem}\label{lem.cov}
For any $s,t \ge 0$ and $x,y \in {\bf R}$,
\[
{\mathbb E}[\psi (s,x) \psi (t,y)]
=
\frac{1}{4 \sqrt{\pi}} 
  \sum_{n \in {\mathbf Z}}
\int_{|s-t|}^{s+t} \sqrt{\frac{1}{l}}  
\exp \Bigl( -\frac{(x-y -n)^2}{4l}  \Bigr) dl.
%\nn
\]
\end{lem}

\Proof
%\QED
%%%%%%%%%%%%%%%%%%%%%%
%
From the explicit expression of $\psi$, we have
\begin{align}
\lefteqn{
{\mathbb E}[\psi (s,x) \psi (t,y)]
}
\nn\\
&=
\int_0^s \int_0^t  drdr'
\int_{S^1}\int_{S^1} dzdz'
\tilde{p}_{s-r} (x-z)  \tilde{p}_{t-r'} (y-z')
{\mathbb E}[\xi (r,z) \xi (r',z')]
\nn\\
&=
\int_0^{s \wedge t} dr \int_{0}^1 \tilde{p}_{s-r} (x-z)  \tilde{p}_{t-r} (y-z) dz
\nn
\nn\\
&=
\int_0^{s \wedge t} dr \sum_{n \in {\mathbf Z}}
\int_{-\infty}^{\infty} p_{s-r} (x-z) p_{t-r} (y+n -z) dz.
\label{eq.covpsi}
\end{align}

A well-known calculation for the heat kernel $p_t$ yields;
\begin{align}
\lefteqn{
\frac{1}{4\pi} 
\int_0^{s \wedge t} dr
\int_{-\infty}^{\infty}
\frac{1}{\sqrt{(s-r)(t-r)} } 
\exp \Bigl( -\frac{(x-z)^2}{4(s-r)}  -\frac{(-z)^2}{4(t-r)} \Bigr)
dz
}
\nn\\
&= 
\frac{1}{2 \sqrt{\pi}}  \int_0^{s \wedge t} dr
(s+t -2r)^{- 1/2} 
\exp \Bigl( -\frac{x^2}{4(s+t-2r)}  \Bigr)
\nn\\
&=
\frac{1}{4 \sqrt{\pi}}  \int_{|s-t|}^{s+t} \sqrt{\frac{1}{l}}  
\exp ( -\frac{x^2}{4l}) dl.
\nn
\end{align}
for all $s, t \ge 0$ and $x \in {\mathbf R}$.
Here, we performed completing the square for the first equality 
and changed variables by $s+t -2r =l$ for the second equality. 
From this, (\ref{eq.covpsi}), and translation invariance,  we prove the lemma.
In the above argument, we used the space-time white noise to calculate the covariance.
However, some readers may think this is not so mathematically rigorous. 
Hence, we will give another proof by using Fourier series.
Remember that both Fourier analysis and Gaussian measure theory work perfectly in $L^2$-setting.

Set $v_0 (x) =1$. For $n >0$,
set also $v_n (x) = \sqrt{2} \cos (2 \pi n x)$ and $v_{-n} (x) = \sqrt{2} \sin (2 \pi n x)$.
Then, $\{ v_n \}_{ n \in {\bf Z}}$ forms an orthonormal basis of $L^2(S^1)$.
For each $n$, $v_n$ is an eigenfunction of $\triangle$ with eigenvalue $4\pi^2 n^2$.
The heat kernel admits Fourier expansion as follows;
\[
\tilde{p}_t (x-y) = \sum_{ n \in {\bf Z}}  e^{- 4\pi^2 n^2 t}  v_n(x) v_n(y).
\]
Using these, we can easily check that
\[
\int_0^1 \tilde{p}_{s-r} (x-z) \tilde{p}_{t-r} (y-z) dz 
=  
\sum_{ n \in {\bf Z}}  e^{- 4\pi^2 n^2 (s-r)}  e^{- 4\pi^2 n^2 (t-r)}  v_n(x) v_n(y).
\]

In a similar way, since $\xi (t,x) dt$  is the increment of an $L^2(S^1)$-cylindrical Brownian motion,
we have
\[
\xi (t,x) dt = \sum_{ n \in {\bf Z}}  db_n (t) v_n (x),
\]
where $(b_n (t))_{t \ge 0}$ are independent copies of the standard real-valued Brownian motion.

If we expand $\psi (t,x) = \sum_{ n \in {\bf Z}}  \hat{\psi}_n (t) v_n (x)$,
then we see the Fourier coefficient $\hat{\psi}_n (t)$ satisfies the following SDE;
\[
d \hat{\psi}_n (t) = 4\pi^2 n^2 \hat{\psi}_n (t) dt + db_n(t), \qquad  \hat{\psi}_n (0)=0.
\]
This SDE has an explicit solution, that is, 
 $\hat{\psi}_n (t) = \int_0^t e^{ -4\pi^2 n^2 (t-r)} db_n (r)$.
Hence, 
\begin{align}
{\mathbb E}[\psi (s,x) \psi (t,y)]
&=
{\mathbb E}  
\Bigl[
\sum_{ n,m \in {\bf Z}} 
\int_0^s e^{ -4\pi^2 n^2 (s-r)} db_n (r)  \cdot \int_0^t e^{ -4\pi^2 m^2 (t-r)} db_m (r)
\cdot v_n(x) v_m (y)
\Bigr]
\nn\\
&=
\sum_{ n \in {\bf Z}} 
\int_0^{s \wedge t}  dr e^{ -4\pi^2 n^2 (s-r)} e^{ -4\pi^2 n^2 (t-r)}  v_n(x) v_n (y)
\nn\\
&=
\int_0^{s \wedge t}  dr
\int_0^1 \tilde{p}_{s-r} (x-z) \tilde{p}_{t-r} (y-z) dz.
\nn
\end{align}
This coincides with (\ref{eq.covpsi}).
Thus, we have obtained the covariance of $\psi$ via Fourier analysis, too. 
\QED

%%%%%%%%%%%%%%%%    Lemma   %%%%%%%%%%

%\vspace{10mm}
%\noindent
%$\spadesuit$~

Next we calculate variance for the two-parameter increment of $\psi$.
Set 
\[
D(s,x; t,y) := {\mathbb E} \Bigl[
\Bigl|
 \psi (t,y)  -  \psi (t,x) - \psi (s,y)  +  \psi (s,x)
     \Bigr|^2
\Bigr].
\]
This quantity plays a very important role in this paper.

%%%%%%%%%%%%%%%%%%%%%%%%%%%%%%%

\begin{lem}\label{lem.est.D^2}
For any $T>0$ and $\kappa \in (0, 1)$, there exists a positive constant 
$c =c(T,\kappa)$ such that 
\[
D(s,x; t,y) \le c |t-s|^{\kappa /2} {\rm dist}_{S^1} (x,y)^{1 -\kappa}
\]
for all $s,t \in [0,T]$ and $x, y \in {\bf R}$.
Here, ${\rm dist}_{S^1} (x,y) = \inf_{n \in {\bf Z}} |x-y+n|$ is the distance on $S^1$. 
(If $|x-y| \le 1/2$, then ${\rm dist}_{S^1} (x,y) =|x-y|$.)
\end{lem}

%%%%%%%%%%%%%%%%%%%%%%%%%%%%%%%%%%%%%%%%%%%%%\vspace{10mm}

\Proof
Due to the periodicity, the invariance
under the translation  and the inversion $x \mapsto -x$, 
we may assume that $s \le t$, $x =0$ and $0<y \le 1/2$
without loss of generality.
In this proof, the positive constant $c$ may change from line to line.

Set $v_t = \psi (t,y)  -  \psi (t,0)$.
Then, we have 
\begin{align}
{\mathbb E} [ v_t v_s]
&=
\frac{1}{2 \sqrt{\pi}} 
  \sum_{n \in {\mathbf Z}}
\int_{|s-t|}^{s+t} \sqrt{\frac{1}{l}} 
\Bigl[
\exp \Bigl( -\frac{n^2}{4l}  \Bigr) - \exp \Bigl( -\frac{(y+n)^2}{4l}  \Bigr) 
\Bigr]
dl.
\nn
\end{align}
Hence, 
\begin{align}
\lefteqn{
D(s,0; t,y)  = {\mathbb E} [ (v_t - v_s )^2]
}
\nn\\
&=
\frac{1}{2 \sqrt{\pi}} 
  \sum_{n \in {\mathbf Z}}
  \Bigl\{ \int_{0}^{2t} -2 \int_{|s-t|}^{s+t} +\int_{0}^{2s}\Bigr\}
   \sqrt{\frac{1}{l}} 
\Bigl[
\exp \Bigl( -\frac{n^2}{4l}  \Bigr) - \exp \Bigl( -\frac{(y+n)^2}{4l}  \Bigr) \Bigr]  dl. 
\nn
\end{align}

%%%%%%%%%%%%%%%%%%%%%%%%%%%%%%%%%%%%%%%%%%%%%\vspace{10mm}

There are two cases. 
{\bf Case (i)}~ the case $t-s \le 2s$.
Then, up to a Lebesgue zero set, 
\[
I_{[0,2t]} -2 I_{[t-s,t+s]} + I_{[0,2s]}
=
2 I_{[0,t-s]} - I_{[2s,t+s]} +I_{[t+s, 2t]},
\]
where $I_A$ stands for the indicator function of $A \subset {\bf R}$.
Note that length of all the intervals on the right hand side is dominated by $t-s$.

{\bf Case (ii)}~ the case $t-s \ge 2s$.
Then, $2t \le 3(t-s)$ and 
\[
|I_{[0,2t]} -2 I_{[t-s,t+s]} + I_{[0,2s]}|
\le
4 I_{[0,3(t-s)]}.
\]

Either way, it is sufficient to prove an inequality of the following form;
\begin{align}
 \sum_{n \in {\mathbf Z}}
  \int_{u }^{ u+ \delta} 
   \sqrt{\frac{1}{l}} 
\Bigl|
\exp \Bigl( -\frac{n^2}{4l}  \Bigr) - \exp \Bigl( -\frac{(y+n)^2}{4l}  \Bigr) \Bigr| dl
\le 
c \delta^{\kappa /2} y^{1 -\kappa}
\label{est.wa.diff}
\end{align}
for all $u \ge 0, \delta \ge 0$ such that $u+\delta \le 3T$.

%%%%%%%%%%%%%%%%%%%%%%% kaki  naosi %%%%%%%%%%%\vspace{10mm}

First, we estimate the $0${th} term in the sum (\ref{est.wa.diff}).
We will change variables from $l$ to $z$ by $y^2 /l =z$ below, then $dl = - (y/z)^2 dz$.
\begin{align}
 \int_{u }^{ u+ \delta} 
   \sqrt{\frac{1}{l}} 
\Bigl(1  - \exp \Bigl( -\frac{y^2}{4l}  \Bigr) \Bigr) dl
&\le 
\int_{0 }^{ \delta} 
   \sqrt{\frac{1}{l}} 
\Bigl(1  - \exp \Bigl( -\frac{y^2}{4l}  \Bigr) \Bigr) dl
\nn\\
&\le 
y \int_{ y^2/\delta}^{\infty}  z^{-3/2}  (1 - e^{-z/4} ) dz.
 \label{eq.zint}
  \end{align}
We will show that
\begin{equation} 
\int_{ r}^{\infty}  z^{-3/2}  (1 - e^{-z/4} ) dz \le c r^{-\kappa/2}
\qquad
\qquad
(0<r<\infty)
 \label{eq.zint2}
\end{equation} 
for some constant $c=c_{\kappa}>0$.
When $r \searrow 0$, there is no problem since the integral is convergent. 
When $r \to \infty$, the integral is dominated by 
$\int_{ r}^{\infty}  z^{-3/2} dz =O (r^{-1/2}) = O (r^{-\kappa/2})$.
This proves (\ref{eq.zint2}).
From (\ref{eq.zint}) and (\ref{eq.zint2}), we have
\begin{equation}\label{est.0th}
\int_{u }^{ u+ \delta} 
   \sqrt{\frac{1}{l}} 
\Bigl(1  - \exp \Bigl( -\frac{y^2}{4l}  \Bigr) \Bigr) dl
\le 
 c  \delta^{\kappa /2} y^{1 -\kappa}.
 \nn
 \end{equation}
Here, the constant $c>0$ does not depend on $y, \delta, u$.

%%%%%%%%%%%%%%%%%%%%%%%%%%%%%%%%%%\vspace{10mm}

Next, we estimate the $\sum_{n >0}$ part of the sum (\ref{est.wa.diff}).
By $(e^{- x^2/4l})^{\prime} = - e^{- x^2/4l} (x/2l)$
and  the mean value theorem,
\[
\exp \Bigl( -\frac{n^2}{4l}  \Bigr) - \exp \Bigl( -\frac{(y+n)^2}{4l}  \Bigr)
\le
 \exp \Bigl( -\frac{n^2}{4l}  \Bigr) \frac{(n+y) y}{2l}
 \le 
 \exp \Bigl( -\frac{n^2}{4l}  \Bigr) \frac{ny}{l},
  \]
where we used $n \ge 1$ and $0 \le y \le 1/2$.

By Schwarz's inequality and change of variables by $n^2/l =z$, we have
\begin{align}
\lefteqn{
 \int_{u }^{ u+ \delta} 
   \sqrt{\frac{1}{l}} 
\Bigl( \exp \Bigl( -\frac{n^2}{4l}  \Bigr)  - \exp \Bigl( -\frac{(n+y)^2}{4l}  \Bigr) \Bigr) dl
}\nn\\
&\le 
yn 
 \int_{u }^{ u+ \delta}  l^{-3/2}
  \exp \Bigl( -\frac{n^2}{4l}  \Bigr)  dl 
 \nn
 \\
  &\le 
  yn 
   \Bigl\{  \int_{u }^{ u+ \delta}   l^{ \kappa -1}   dl \Bigr\}^{\frac12} 
              \Bigl\{  \int_{u }^{ u+ \delta}  l^{-(\kappa +2)}  
           \exp \Bigl( -\frac{n^2}{2l}  \Bigr)   dl \Bigr\}^{\frac12}     
   \nn\\
  &\le  
  cy \delta^{\kappa /2} n 
   \Bigl\{  \int_{u }^{ u+ \delta}  l^{-(\kappa +2)}  
           \exp \Bigl( -\frac{n^2}{2l}  \Bigr)   dl \Bigr\}^{\frac12} 
              \nn\\
               &\le 
                 cy \delta^{\kappa/2} n \cdot n^{-(1+\kappa)}  
                  \Bigl\{
                    \int_{n^2/ 3T}^{\infty}  z^{\kappa} e^{-z/2}dz    \Bigr\}^{\frac12}
                    \nn\\
                      &\le 
                          cy \delta^{\kappa /2}  \exp \Bigl( -\frac{n^2}{24T}  \Bigr),
                  %
                  %\label{int.ldelta}
 \nn                   
  \end{align}
since there exists $c>0$ such that $z^{\kappa} e^{-z/2} \le c e^{-z/4}$ for all $z \ge 0$.
Since
\[
\sum_{n=1}^{\infty} \exp \Bigl( -\frac{n^2}{24T}  \Bigr)
\le
\int_0^{\infty} \exp \Bigl( -\frac{x^2}{ 24T} \Bigr) dx \le \sqrt{6 \pi T},
\]
we  have
\[
\sum_{n=1}^{\infty}
\int_{u }^{ u+ \delta} 
   \sqrt{\frac{1}{l}} 
\Bigl( \exp \Bigl( -\frac{n^2}{4l}  \Bigr)  - \exp \Bigl( -\frac{(n+y)^2}{4l}  \Bigr) \Bigr) dl
\le 
 cy \delta^{\kappa /2}. %
\]
Using the condition that $0 \le y \le 1/2$,
we can also estimate the $\sum_{n < 0}$ part of (\ref{est.wa.diff}) with a slight modification.
Thus we have shown (\ref{est.wa.diff}) and consequently Lemma \ref{lem.est.D^2}.
\QED

%%%%%%%%%%%%%%%%    Lemma   %%%%%%%%%%

%\vspace{10mm}
%\noindent
%$\spadesuit$~

By the previous lemma and
Kolmogorov-\v{C}encov's continuity criterion,
we can obtain the regularity of sample sheet of $\psi$.
\begin{cor}\label{cor.kolm}
{\rm (i)}~ For any $T>0$, there exists a positive constant $c=c_T$
such that
\begin{eqnarray}
{\mathbb E} \bigl[
\bigl|
   \psi (t,x) -  \psi (s,x)
     \bigr|^2
\bigr]
\le c |t-s|^{1/2},
\qquad
{\mathbb E} \bigl[
\bigl|
   \psi (t,x) -  \psi (t,y)
     \bigr|^2
\bigr]
\le c |y-x|
\nn
\end{eqnarray}
holds for all $s,t \in [0,T]$ and $x,y \in [0,1]$.
\\
\noindent 
{\rm (ii)}~ For any $\beta <1/4$ and $\beta' <1/2$,
 $(t,x) \mapsto \psi(t,x)$ is $(\beta, \beta')$-H\"older continuous almost surely.
 That is, 
 \[
 \| \psi\|_{H;(\b,\b')}  :=
 \sup \Bigl\{ \frac{| \psi(s,x) - \psi(t,y)|}{|t-s|^{\beta}+ |x-y|^{\beta'}}  
 ~|~  (s,x)\neq (t,y) 
 \Bigr\} 
  < \infty  \qquad \mbox{a.s.}
 \]
 Here, $\sup$ runs over all $(s,x), (t,y) \in  [0,T]\times [0,1]$ such that $ (s,x)\neq (t,y)$.
\end{cor}

%%
%%%%%

\Proof
By repeating a similar calculation as in the proof of Lemma \ref{lem.est.D^2},
we can prove the first assertion. (This is actually  easier than  Lemma \ref{lem.est.D^2}.
Note that one should not take limit ($\kappa \searrow 0$ or $\kappa \nearrow 1$) in 
Lemma \ref{lem.est.D^2}, since the constant $c$ may depend on $\kappa$.)

To prove the second assertion, 
we use Kolmogorov-\v{C}encov's continuity criterion (see e.g. Theorem 1.4.4, p. 36, Kunita \cite{ku}).
From the first assertion and Lemma \ref{lem.est.D^2}, 
we see that,
for any fixed $\kappa \in (0,1) \cap {\bf Q}$ and $\gamma \in {\bf N}$,
\begin{eqnarray*}
{\mathbb E} \bigl[
\bigl|
   \psi (t,x) -  \psi (s,x)
     \bigr|^{\gamma}
\bigr]
\le c |t-s|^{\kappa \gamma/4},
\qquad
{\mathbb E} \bigl[
\bigl|
   \psi (t,x) -  \psi (t,y)
     \bigr|^{\gamma}
\bigr]
\le c |y-x|^{(1-\kappa) \gamma /2},
\nn
\\
 {\mathbb E} \bigl[
\bigl|
 \psi (t,y)  -  \psi (t,x) - \psi (s,y)  +  \psi (s,x)
     \bigr|^{\gamma}
\bigr]
\le
c |t-s|^{\kappa \gamma/4} |x-y|^{(1 -\kappa)\gamma /2}.
\end{eqnarray*}
Here, $c>0$ may depend on $\gamma$.
Hence, by the continuity criterion,
there exists a modification of $\psi$ (which is again denoted by the same symbol)
which is a.s. $(\dl, \dl')$-H\"older continuous for any $(\dl, \dl')$ such that
\[
\dl < \frac{\kappa}{4} - \frac{1}{\gamma}, \qquad 
\dl' < \frac{1- \kappa}{2} - \frac{1}{\gamma}.
\]
First choose $\kappa_1$ and $\gamma_1$ so that $\b < \kappa_1/4 - 1/\gamma_1$.
We set $\dl_1 =\b$ and 
take $\dl_1'$ so that $0<\dl_1' < (1- \kappa_1)/2 - 1/\gamma_1$.
Next
choose $\kappa_2$ and $\gamma_2$ so that $\b'  < (1- \kappa_2)/2 - 1/\gamma_2$.
We set $\dl_2' =\b'$ and take $\dl_2$ so that $0<\dl_2 < \kappa_2/4 - 1/\gamma_2$.
We now have two modifications, but they coincide a.s. anyway.
Hence, we have
\begin{align}
|\psi (s,x)  -  \psi (t,y) | &\le
 |\psi (s,x)  -  \psi (t,x) | + |\psi (t,x)  -  \psi (t,y) | 
\nn\\
&\le 
\|\psi\|_{H;(\dl_1,\dl_1')} |t-s|^{\b} + \|\psi\|_{H;(\dl_2,\dl_2')} |x-y|^{\b'},
 \qquad \mbox{a.s.}
 \nn
\end{align}
This proves the second assertion.
\QED

%%%%%%%%%%%%%%%%    Lemma   %%%%%%%%%%

%
%
Now we check that $\psi$ satisfies a condition of Coutin-Qian type 
(see Definition 4.4.1, Lyons and Qian \cite{lq}).
Ours is slightly weaker than the one in Definition 4.4.1, \cite{lq},
but practically there is no problem.
This kind of condition implies 
that 
the lift of
the dyadic piecewise linear approximation of Gaussian process 
converges in the geometric rough path space.

%%%%%%%

\begin{lem}\label{lem.cq1}
For any $T>0$, there exists a positive constant $c=c_T$ such that
\begin{equation}\label{est.cq1}
\bigl|
{\mathbb E}[ \{\psi(t, x+h) -\psi(t, x)\}   \{\psi(s, y+h) -\psi(s, y)\} ] 
\bigr|
\le
 \frac{c h^2}{y-x} 
\end{equation}
holds for all $s,t \in [0,T]$ and all $0 \le x<y \le 1$, $h>0$ such that $2h \le y-x \le 1/2$.
\end{lem}

%%%%%%%%

\Proof 
In this proof, the constant $c>0$ may change from line to line.
As before, it is sufficient to prove the lemma when $x=0$.

Let $h>0$ and let
$g$ be a real-valued, $C^2$-function defined on a certain interval
which includes $[y-h, y+h]$.
By Taylor's theorem, 
\begin{align}
|g(y+h) + g(y-h) -2g(y) |
&=
h^2 
\Bigl|
\int_0^1  (1 - \theta) \{ g^{\prime\prime} (y +\theta h) + g^{\prime\prime} (y -\theta h)\} d\theta
\Bigr|
\nn\\
&=
h^2 
\Bigl|
\int_{-1}^1  \{(1 - \theta) \wedge (1 + \theta) \}
g^{\prime\prime} (y +\theta h)  d\theta
\Bigr|
\nn\\
&\le
h^2   \sup_{y-h \le \eta \le y+h} | g^{\prime\prime} (\eta ) |.
%\int_{-1}^1
%
%| g^{\prime\prime} (y +\theta h) | d\theta
%
\label{eq.tay.2nd}
\end{align}
We will use (\ref{eq.tay.2nd}) for $g (y) =\exp (- y^2/ 4l )$.
It is easy to see that
$g^{\prime} (y) = - (y/2l) \exp (- y^2/ 4l )$
and 
$g^{\prime\prime} (y) = \{ (y^2 /4l^2) - (1/2l)\} \exp (- y^2/ 4l )$.

From the covariance formula for $\psi$ in Lemma \ref{lem.cov},
the left hand side of (\ref{est.cq1}) is dominated by a constant multiple of 
\begin{eqnarray}
  \sum_{n \in {\mathbf Z}}
\int_{0}^{2T} \sqrt{\frac{1}{l}}  
\Bigl|
\exp \Bigl( -\frac{(y+h +n)^2}{4l}  \Bigr) + \exp \Bigl( -\frac{(y -h +n)^2}{4l}  \Bigr) 
-2\exp \Bigl( -\frac{(y+n)^2}{4l}  \Bigr) 
\Bigr| dl
\nn
%
%\\
%=:\sum_{n \in {\mathbf Z}}A(n).
%\label{est.cq1.2}
\end{eqnarray}
We will denote by $A(n)$ the $n$th summand in the above sum.

First we estimate $A(0)$.
From (\ref{eq.tay.2nd}) and $0 \le h \le y/2$, we see that
\begin{align}
A(0) 
&\le
h^2 \int_{0}^{2T} \sqrt{\frac{1}{l}}  
\Bigl\{  \frac{ (3y/2)^2}{4 l^2} + \frac{1}{2l} \Bigr\}
\exp \Bigl( -\frac{(y/2)^2}{4l}  \Bigr) 
dl
\nn\\
&\le
ch^2 \int_{0}^{2T} \sqrt{\frac{1}{l}}  
\Bigl( \frac{ y^2}{ l^2} + \frac{1}{l} \Bigr)
e^{ -y^2 /16l}   
dl.
\nn
\end{align}
Change variables from $l$ to $z$ by $y^2/l =z$.
Then, $dl = - (y^2 /z^2) dz$ and we have
\begin{align}
A(0) 
&\le
\frac{ch^2}{y}
\int_{y^2 /2T}^{\infty}  (z^{1/2} + z^{- 1/2}) e^{-z /16}dz \le \frac{ch^2}{y}.
\label{est.A0}
\end{align}
Note that the integral above is convergent on $(0, \infty)$.

Next we estimate $A(n)$ for $n>0$.
Using (\ref{eq.tay.2nd}) with $y+n$ instead of $y$ and repeating a similar computation, we have
\begin{align}
\sum_{n=1}^{\infty}
A(n) 
&\le
h^2
\sum_{n=1}^{\infty}
\int_{0}^{2T} \sqrt{\frac{1}{l}}  
\Bigl\{  \frac{ (3y/2 +n )^2}{4 l^2} + \frac{1}{2l} \Bigr\}
\exp \Bigl( -\frac{(y/2 +n)^2}{4l}  \Bigr) 
dl
\nn\\
&\le
c h^2 \sum_{n=1}^{\infty}
\int_{0}^{2T} \sqrt{\frac{1}{l}}  
\Bigl( \frac{ n^2}{ l^2} + \frac{1}{l} \Bigr)
e^{ - n^2 /4l}   
dl
\nn\\
&\le
\sum_{n=1}^{\infty}
\frac{c h^2}{n}
\int_{n^2 /2T}^{\infty}  (z^{1/2} + z^{- 1/2}) e^{-z /4}dz 
\nn\\
&\le 
\sum_{n=1}^{\infty}
\frac{c h^2}{n}
\int_{n^2 /2T}^{\infty}  z^{- 1/2} e^{-z /8}dz 
\nn\\
&\le 
\sum_{n=1}^{\infty}
\frac{ch^2}{n^2}
\int_{0}^{\infty}  e^{-z /8}dz 
\le ch^2 \le \frac{ch^2}{y}.
\label{est.An>0}
\end{align}
Here, we used $z+1 \le ce^{z/8} ~(z>0)$ for some $c>0$ and $1 \le (2y)^{-1}$.

Using the assumption that $y \le 1/2$, we can also prove 
$\sum_{n=1}^{\infty} A(-n) < ch^2 /y$ essentially in the same way.
Combining this with (\ref{est.A0}) and  (\ref{est.An>0}), we have shown the lemma.
\QED

%%%%%%%%%%%%%%%%%%%%%%%%%%%%%%%%%%%%%%%%%%%%%%%%%%%%%%%%
%%%%%%  Coutin-Qian like condition with time difference
%%%%%%%%%%%%%%%%%%%%%%%%%%%%%%%%%%%%%%%

%%%%%%%%%%%%%%%    Lemma   %%%%%%%%%%

%\vspace{10mm}
%\noindent
%$\spadesuit$~

The following is a generalized version of Coutin-Qian's condition 
in the sense that the regularity with respect to $t$-variable is also taken into account.
This is a key technical lemma and 
 will play a crucial role in the next section when we prove convergence of 
dyadic polygonal approximation on the path space over the geometric rough path space.

%%%%%%%

\begin{lem}\label{lem.cq2}
For any $T>0$ and $\kappa \in (0,1)$, 
there exists a positive constant $c=c(T,\kappa)$ such that
\begin{align}
\lefteqn{
\Bigl|
{\mathbb E}
\Bigl[ 
\bigl\{ (\psi(t, x+h) -\psi(t, x) ) -   (\psi(s, x+h) -\psi(s, x) )\bigr\} 
}\nn\\
& \qquad
\times
\bigl\{ (\psi(t, y+h) -\psi(t, y) ) -   (\psi(s, y+h) -\psi(s, y) )\bigr\}
\Bigr] 
\Bigr|
&\le
 \frac{c |t-s|^{\kappa /2} h^2}{(y-x)^{1 +\kappa}} 
\label{est.cq2}
\end{align}
holds for all $s,t \in [0,T]$ and all 
$0 \le x<y \le 1$, $h>0$ such that $2h \le y-x \le 1/2$.
\end{lem}

%\vspace{5mm}

\Proof
In this proof, the constant $c>0$ may change from line to line.
As before, it is sufficient to prove the lemma when $x=0$.

From the covariance formula in Lemma \ref{lem.cov},
the left hand side of (\ref{est.cq1}) is dominated by a constant multiple of 
$\sum_{n \in {\bf Z}} |B(n)|$, where
\begin{align}
B(n) 
&=
 \Bigl\{ \int_{0}^{2t} -2 \int_{|s-t|}^{s+t} +\int_{0}^{2s}\Bigr\}
 \sqrt{\frac{1}{l}}  
\nn\\
& \times
\Bigl[
\exp \Bigl( -\frac{(y+h +n)^2}{4l}  \Bigr) + \exp \Bigl( -\frac{(y -h +n)^2}{4l}  \Bigr) 
-2\exp \Bigl( -\frac{(y+n)^2}{4l}  \Bigr) 
\Bigr] dl.
\label{est.cq2.1}
\end{align}
It is obvious from (\ref{est.cq2.1}) that we may assume $s \le t$ without loss of generality.

For the same reason as in the proof of Lemma \ref{lem.est.D^2} (see {\bf Cases (i)(ii)}),
it is sufficient to prove the following inequality:
There exists $c>0$ such that 
%
%for any $u \ge 0$, $\delta \ge 0$ with $u +\delta \le 3T$, it holds that
%
%
\begin{equation}
\sum_{n \in {\bf Z}} C(n) \le  \frac{c \delta^{\kappa /2} h^2}{y^{1 +\kappa}}
\qquad
\qquad
(u \ge 0, \delta \ge 0 \mbox{ with } u +\delta \le 3T),
\label{est.cq2.2}
\end{equation}
where
\begin{align}
C(n) 
=
 \int_{u}^{u +\delta}
 \sqrt{\frac{1}{l}}  
\Bigl|
\exp \Bigl( -\frac{(y+h +n)^2}{4l}  \Bigr) + \exp \Bigl( -\frac{(y -h +n)^2}{4l}  \Bigr) 
-2\exp \Bigl( -\frac{(y+n)^2}{4l}  \Bigr) 
\Bigr| dl.
\nn
%
%\label{est.cq2.3}
\end{align}
%

%%%%%%%%%%%%%%%%%%%%%%%%%%%]

Using (\ref{eq.tay.2nd}), we calculate $C(0)$
in the same way as for $A(0)$ in  the proof of Lemma \ref{lem.cq1}.
\begin{align}
C(0) 
&\le
h^2 \int_{u}^{u +\delta} 
l^{-1/2} 
\Bigl\{  \frac{ (3y/2)^2}{4 l^2} + \frac{1}{2l} \Bigr\}
\exp \Bigl( -\frac{(y/2)^2}{4l}  \Bigr) 
dl
\nn\\
&\le
ch^2 \int_{u}^{u+\delta} 
l^{-3/2} 
\Bigl( \frac{ y^2}{ l} + 1 \Bigr)
e^{ -y^2 /16l}   
dl
\nn\\
&\le
ch^2
\Bigl\{ \int_{u }^{ u+ \delta}  (l^{\frac{\kappa -1}{2}})^2 dl\Bigr\}^{\frac12} 
      \Bigl\{  \int_{u }^{ u+ \delta}  (l^{- \frac{\kappa +2}{2}})^2   
      \Bigl( \frac{ y^2}{ l} + 1 \Bigr)^2     
          e^{ -y^2 /8l}  dl \Bigr\}^{\frac12} 
        \nn
 %\\
%          \\
           &\le  
             ch^2 \delta^{\kappa /2}
     \Bigl\{ \frac{1}{y^{2+2\kappa}}  \int_{y^2/3T }^{ \infty}  z^{\kappa} (z+1)^2  
       e^{ -z /8}  dl \Bigr\}^{\frac12} 
         \le 
          \frac{ch^2 \delta^{\kappa /2}}{y^{1+\kappa}}.
        \nn
\end{align}
In the last line above, we changed variables by $y^2/l =z$ again.

Let us estimate $C(n)$ for $n>0$.
In the same way as in the proof of Lemma \ref{lem.cq1},
\begin{align}
\sum_{n=1}^{\infty}
C(n) 
&\le
h^2
\sum_{n=1}^{\infty}
\int_{u}^{u+\delta} \sqrt{\frac{1}{l}}  
\Bigl\{  \frac{ (3y/2 +n )^2}{4 l^2} + \frac{1}{2l} \Bigr\}
\exp \Bigl( -\frac{(y/2 +n)^2}{4l}  \Bigr) 
dl
\nn\\
&\le
c h^2
\sum_{n=1}^{\infty}
\int_{u}^{u+\delta} l^{-3/2} 
\Bigl( \frac{ n^2}{ l} + 1 \Bigr)
e^{ - n^2 /4l}   
dl
\nn\\
&\le
c h^2
\sum_{n=1}^{\infty}
\Bigl\{ \int_{u }^{ u+ \delta}  (l^{\frac{\kappa -1}{2}})^2 dl\Bigr\}^{\frac12} 
      \Bigl\{  \int_{u }^{ u+ \delta}  (l^{- \frac{\kappa +2}{2}})^2   
      \Bigl( \frac{ n^2}{ l} + 1 \Bigr)^2     
          e^{ -n^2 /4l}  dl \Bigr\}^{\frac12} 
        \nn
 \\
  &\le c h^2
    \delta^{\kappa /2} \sum_{n=1}^{\infty}
  \Bigl\{ \frac{1}{n^{2+2\kappa}}  \int_{n^2/3T }^{ \infty}  z^{\kappa} (z+1)^2  
       e^{ -z /4}  dl \Bigr\}^{\frac12}
        \nn\\
        &\le
           c h^2
           \delta^{\kappa /2} \sum_{n=1}^{\infty} \frac{1}{n^{1+\kappa}} 
           \le c h^2\delta^{\kappa /2}
           \le \frac{c h^2\delta^{\kappa /2}}{ y^{1 +\kappa}}.   
           \nn               
           \end{align}
Using the condition $y \le 1/2$, we can also prove the same estimate for $\sum_{n=1}^{\infty}C(-n)$.
Thus, we have shown (\ref{est.cq2.2}).
\QED

%\newpage
%%%%%%%%%%%%%%%%%%%%%%%%%%%%%%%%%%%%%%%%%%%%%%%%%%%%%%%%%%%%%%%%%%%%%%%%%%%%%%%%
%%%%%%%%%%%%%        Section     %%%%%%%%%%%%%%%%%%%%%%%%%%%%%%%%%%%%%%
%%%%%%%%%%%%%%%%%%%%%%%%%%%%%%%%%%%%%%%%%%%%%%%%%%%%%%%%%%%%%%%%%%%%%%%%%%%%%%%%

\section{Dyadic polygonal approximation}
\label{sec.dyconv}

%$\spadesuit$~
In this section we again consider the multi-dimensional case of (\ref{eq.psi.spde})
; $\psi =(\psi^1, \ldots, \psi^d)$.
Here,  $\psi^1, \ldots, \psi^d$ are independent copies of 
the one-dimensional process studied in the previous section.
For the dyadic partition ${\cal P}_k = \{ i /2^k ~|~ 0 \le i \le 2^k\}$ 
of $[0,1]$ with $k=1,2,\ldots$,
we denote by $\psi (k) (t,x)$ the piecewise linear approximation in $x$
associated with ${\cal P}_k$ for each fixed $t$.
That is, for all $t$ and $i$,
$\psi (k) (t, i2^{-k}) =\psi(t, i2^{-k})$
and $x \mapsto \psi (k) (t,x)$ is linear on each $[(i-1)2^{-k}, i2^{-k}]$.
The main objective of this section is to prove that the spatial lift 
$\Psi(k)= {\cal L}_2(\psi(k))$ of $\psi (k)$
converges as $k \to \infty$ in ${\cal P}^B_{\b,m} G\Omega^B_{\a, m} ({\bf R}^d)$
for suitable parameters $\b, \a, m$.

%%%%%%%%%%%%%%%%%%%%%%  equivalence of chaos %%%%%%%%%%%%%%%%%%%%%%%%%%%

%\vspace{10mm}
%\noindent
%$\spadesuit$~

Before we calculate the dyadic polygonal approximations, 
we provide a proposition on equivalence of a Banach space-valued Wiener chaos.
This proposition is quite useful for our purpose and will be used frequently.
It is well-known that, on a fixed real-valued Wiener chaos, all $L^p$-norm are equivalent ($1<p<\infty$).
This is still true in the case of a Banach space-valued Wiener functionals. 
(For instance, see Friz and Victoir \cite{fv07} or Maas \cite{maas})
The following is a quantified version of this fact (Lemmas 2 and 3, \cite{fv07}). 
It is worth noting that it holds for any real Banach space ${\cal V}$
without any additional condition.

\begin{prop}\label{pr.equi.chao}
Let $({\cal X},{\cal H},\mu)$ be an abstract Wiener space. 
For  a real Banach space ${\cal V}$,
denote by ${\cal C}_n ({\cal V})$  the ${\cal V}$-valued $n$th inhomogeneous Wiener chaos.
\\
\noindent
{\rm (i)}~
Restricted on a fixed Wiener chaos ${\cal C}_n ({\cal V})$, 
all $L^p$-norm are equivalent ($1<p<\infty$).
\\
\noindent
{\rm (ii)}~
For any $n$, there exists a positive constant $c=c_n$ such that 
\[
\|  Z\|_{ L^p  }  \le  \|  Z\|_{ L^q }  \le c_n  (q-1)^{n/2} \|  Z\|_{ L^p }  
\]
holds for all $2 \le p \le q <\infty$ and $Z \in {\cal C}_n ({\cal V})$.
\end{prop}

Let 
$g \in C_0 ([0,T], {\cal X})$, where ${\cal X} = C^{B; \a,m} ([0,1], {\bf R}^d)$.
As before, we write 
$G^1(t;x,y) =G^1_t(x,y) =g(t,y) -g(t,x)$.
By slightly abusing notations, we set 
\begin{align}
\| g\|_{B;\b,\a,m}^m :&=  \| G^1\|_{B;\b,\a,m}^m
\nn\\
&=
\iint_{ {\cal S} (T)}  dsdt 
\iint_{  {\cal S} }
\frac{ |g(t,x) -g(s,x) - g(t,y) +g(s,y) |^{m} }{|t-s|^{1 +\b m} |x-y|^{1 +\alpha m} }
dxdy.
\nn
\end{align}
Recall that we assumed $g(0,x) \equiv 0$.

Assume that $\a \in (1/3, 1/2)$ and $\b >0$ with $4\b <1-2\a$.
Then, we can find $\kappa \in (0,1)$ such that $\a <(1 -\kappa)/2$ and $\b <\kappa /4$ hold.
Assume further that $m \ge 1, ~\b >1/m$, and $a -1/m >1/3$.

By Corollary \ref{cor.kolm}, we can easily estimate $(\b,m)$-Besov norm of $t \mapsto \psi(t,0)$.
It holds that ${\mathbb E} [\| \psi (\,\cdot\,,0) \|_{B;\b,m}^m ] <\infty$
and hence $\| \psi (\,\cdot\,,0) \|_{B;\b,m} <\infty$, a.s.
Note that $\psi(t,0) =\psi(k)(t,0)$ for all $t$ and $k$.

By Lemma \ref{lem.est.D^2}, we have 
${\mathbb E} [\| \psi \|_{B;\b,\a,m}^m ] <\infty$ and hence 
$\| \psi \|_{B;\b,\a,m} <\infty$, a.s.
In a similar way,
${\mathbb E} [\| \psi (k) \|_{B;\b,\a,m}^m ] <\infty$ for each fixed $k$.
As a result, the laws of $\psi$ and $\psi (k)$
are Gaussian measures on 
$C_0^{B; \b,m} ([0,T], {\cal X})$ with  ${\cal X} = C^{B; \a,m} ([0,1], {\bf R}^d)$.

Now we prove the convergence of $\{\psi(k) \}_{k=1}^{\infty}$ as $k \to \infty$
in $(\b,\a,m)$-Besov norm.

\begin{lem}\label{lem.Bcon1}
Assume
$\a \in (1/3, 1/2)$ and $0< 4 \b < 1- 2\a$.
In addition, assume $m > 1$ satisfies that $\a -(1/m) >1/3$ and $\b >1/m$.
Then, there are constants  $c>0$, $\eta \in (0,1)$ independent of $k$ such that 
\begin{equation} \label{eq.Bcon1}
 {\mathbb E} \Bigl[  \| \psi (k+1) - \psi(k) \|_{B;\b,\a,m}^m  \Bigr]^{1/m}  \le c \eta^k
\end{equation}
for all $k$.
In particular, $\{\psi(k) \}_{k=1}^{\infty}$ converges to $\psi$ 
in $(\b,\a,m)$-Besov topology almost surely 
and  in $L^p$ for all $p \in (1, \infty)$.
\end{lem}

\Proof 
First, we give an integration formula for later use. 
This is useful when we estimate Besov norms.
For any $T>0$, there exists a constant $c_T >0$ such that  
\begin{equation}\label{eas.int.form}
\iint_{{\cal S}} \frac{ \ve \wedge |x-y|^a}{|x-y|^b} dxdy 
\le 
\frac{ \ve^{(a-b+1)/a} } {(a-b+1)(b-1)},
\qquad
\quad
(0 \le \ve \le 1,~ a>b-1>0).
\end{equation}
To check this formula, change variables by $u=s,~v=t-s$.
Then, the integral domain becomes $\{0<u<1,~ 0<v<1, ~ u+v <1 \}$.
The rest is easy.

In this proof the positive constant $c$ may change from line to line.
For $k=1,2,\ldots$ and $1 \le j \le 2^k$,  we set $\D_{j}^{k} \psi_t = \psi (t, j/2^k) - \psi (t, (j-1)/2^k) $.
We write $\lm (k) = \psi(k+1)-\psi(k)$ for simplicity.

When $(j-1)/2^k \le x \le j/2^k$, 
\begin{align}
\lm (k) (t,x) 
&= 
2^k \{ (x - \frac{j-1}{2^k}) \wedge (\frac{j}{2^k} - x) \} 
\{
2 \psi (t, \frac{2j-1}{2^{k+1}}) - \psi (t, \frac{j-1}{2^k}) - \psi (t,  \frac{j}{2^k})
\}
\nn\\
&=
2^k \{ (x - \frac{j-1}{2^k}) \wedge (\frac{j}{2^k} - x) \} 
 ( \D_{2j-1}^{k+1} \psi_t - \D_{2j}^{k+1}\psi_t ).
 \nn
\end{align}

\noindent
This is just a product of  functions in $t$ and  in $x$. 
Hence, when $(j-1)/2^k \le x,y \le j/2^k$, we have
\begin{align}
\lefteqn{
\lm (k) (t,y) -\lm (k) (t,x) -\lm (k) (s,y) +\lm (k) (s,x) 
}
\nn\\
&=
2^k \{ (y - \frac{j-1}{2^k}) \wedge (\frac{j}{2^k} - y) - (x - \frac{j-1}{2^k}) \wedge (\frac{j}{2^k} - x) \} 
\nn\\
& \qquad 
\times
\{   ( \D_{2j-1}^{k+1} \psi_t - \D_{2j-1}^{k+1}\psi_s ) 
  - ( \D_{2j}^{k+1} \psi_t - \D_{2j}^{k+1}\psi_s ) \}.
\nn
\end{align}

\noindent
By using Lemma \ref{lem.est.D^2}, we can estimate the variance of this Gaussian random variable as follows: 
\begin{align}
\lefteqn{
{\mathbb E} [ |\lm (k) (t,y) -\lm (k) (t,x) -\lm (k) (s,y) +\lm (k) (s,x) |^2]
}
\nn\\
&\le
c 2^{2k} |x-y|^2  (2^{-k})^{1-\kappa} |t-s|^{\kappa /2} 
\le 
c |t-s|^{\kappa /2}  |x-y|^{1 -\kappa}.
\nn
\end{align}

Next we consider the case that $x$ and $y$ are in distinct subintervals.
We may assume there exist $j<l$ such that $(j-1)/2^k \le x  \le j/2^k \le (l-1)/2^k \le y \le l/2^k$.
Since $\lm (k) (t,i/2^k) =0$ for all $i$,
\begin{align}
\lefteqn{
\lm (k) (t,y) -\lm (k) (t,x) -\lm (k) (s,y) +\lm (k) (s,x) 
}
\nn\\
&=
\{  \lm (k) (t,y) -\lm (k) (t, (l-1)/2^k) -\lm (k) (s,y) +\lm (k) (s, (l-1)/2^k) \}
\nn\\
&
\quad + \{ \lm (k) (t,  j/2^k) -\lm (k) (t,x) -\lm (k) (s,  j/2^k) +\lm (k) (s,x)  \}.
\nn
\end{align}
Now we can use the result for the previous case to obtain 
\begin{align}
{\mathbb E} [ |\lm (k) (t,y) -\lm (k) (t,x) -\lm (k) (s,y) +\lm (k) (s,x) |^2]
\le
c |t-s|^{\kappa /2}  (2^{-k})^{1 -\kappa}.
\nn
\end{align}
From these we can easily see that, for all $s,t$ and $x,y$, 
\begin{eqnarray}
\lefteqn{
{\mathbb E} [ |\lm (k) (t,y) -\lm (k) (t,x) -\lm (k) (s,y) +\lm (k) (s,x) |^2]
}
\nn\\
&\le&
{\mathbb E} 
\Bigl[  
|\Psi(k+1)^1_t (x,y)  - \Psi(k)^1_t (x,y)  - \Psi(k+1)^1_s (x,y)  +\Psi(k)^1_s (x,y) |^2
\Bigr]
\nn\\
&\le&
c |t-s|^{\kappa /2} \bigl\{ (2^{-k})^{1 -\kappa}  \wedge    |x-y|^{1 -\kappa} \bigr\}.
\label{est.lmd}
\end{eqnarray}

Using (\ref{eas.int.form}) and choosing $\kappa \in (0,1)$ 
so that $\b <\kappa /4$ and $\a < (1 - \kappa)/2$, 
we have 
\begin{align}
\lefteqn{
{\mathbb E} \Bigl[  \| \psi (k+1) - \psi(k) \|_{B;\b,\a,m}^m  \Bigr]
}
\nn\\
& \le 
c \int_0^T \int_0^t |t-s|^{ (\kappa m /4) -(1+\b m)} dsdt
\int_0^1 \int_0^y  \frac{  (2^{-k})^{ (1 -\kappa) m/2} \wedge   |x-y|^{ (1 -\kappa) m/2} }{  |x-y|^{1 +\a m}}   dxdy
\nn\\
& \le 
c 
\Bigl( \frac12 \Bigr)^{ (1 -\kappa -2\a) mk/2}.
\nn
\end{align}
Thus, we have shown (\ref{eq.Bcon1}) with $\eta = 2^{ - (1 -\kappa -2\a) /2} \in (0,1)$.
Almost sure convergence and $L^m$-convergence are immediate from (\ref{eq.Bcon1}).
By Proposition \ref{pr.equi.chao}, convergence is also in $L^p$ for all $p~(1<p<\infty)$.
\QED

%%%%%%%%%%%%%%%%%%%%%%%%%%%%%%%%%%%%%%%%%%%%%%%%%%%%%%%%%%%%%%%%%%%%%%%%%%%%%%%%%%%%
%%%      Second level path
%%%%%%%%%%%%%%%%%%%%%%%%%%%%%%%%%%%%%%%%%%%%%%%%%%%%%%%%%%%%%%%%%%%%%%%%%%%%%%%%%%%%%%

%

As usual we define the second level path of $\psi(k)$ as follows.
For all $0 \le x \le y \le 1$, we set
$$
\Psi (k)^2(t;x,y) := \int_x^y  \{ \psi(k)(t, u) -\psi(k)(t, x) \} \otimes d_u \psi(k)(t, u).
$$
We will often write $\Psi (k)^2_t (x,y)$ for the left hand side when there is no possibility of confusion.
By a slight abuse of notation, 
${\cal L}_2$ also stands for this natural spatial lift, i.e.,
${\cal L}_2 ( \psi(k)) = \Psi (k) =(\psi(k), \Psi(k)^2) =(\psi( \,\cdot\, , 0),\Psi(k)^1, \Psi(k)^2)$.

For an appropriate choice of the parameters,  $\Psi (k)^2$ is a random variable taking its values in 
$C_0^{B; 2\b,m/2} ([0,T], {\cal X})$.
Here,  ${\cal X}$ is the closure of $C^{B; 2\a,m/2}_0 ({\cal S}, {\bf R}^d \otimes {\bf R}^d)$
with respect to $(2\a,m/2)$-Besov norm.
Recall that, for $G \in C_0^{B; 2\b,m/2} ([0,T], {\cal X})$,
its $(2\b, 2\a, m/2)$-Besov norm is given by
\begin{align}
\| G\|_{B;2\b,2\a,m/2}^{m/2}
:&= 
\iint_{ {\cal S} (T)}  dsdt 
\iint_{  {\cal S} }
\frac{ |G(t;x,y) -G(s;x,y)  |^{m/2} }{|t-s|^{1 +\b m } |x-y|^{1 +\alpha m} }
dxdy.
\nn
\end{align}
%

%%%%%%%%%%   Lemma  lem.Bcon2 %%%%%%%%%%%%%

%
Our purpose in this section is to prove that $\{ \Psi (k)^2\}^{\infty}_{k=1}$ converges as $k \to \infty$
in $(2\b, 2\a,m/2)$-Besov norm for suitable parameters $m, \a, \b$.

\begin{lem}\label{lem.Bcon2}
Assume
$\a \in (1/3, 1/2)$ and $0< 4 \b < 1- 2\a$.
In addition, assume $m > 2$ satisfies that $\a -(1/m) >1/3$ and $\b >1/m$.
Then, there are constants  $c>0$, $\eta \in (0,1)$ independent of $k$ such that 
\begin{equation}\label{ineq.Bcon2}
 {\mathbb E} \Bigl[  \| \Psi (k+1)^2 - \Psi(k)^2 \|_{B;2\b,2\a,m/2}^{m/2}  \Bigr]^{2/m}  \le c \eta^k
\end{equation}
for all $k$.
In particular, $\{\Psi(k)^2 \}_{k=1}^{\infty}$ converges 
in $(2\b, 2\a,m/2)$-Besov norm almost surely 
and  in $L^{p}$ for all $p \in (1,\infty)$.
Consequently,  
under the assumptions of the lemma,
$\{ \Psi (k)\}^{\infty}_{k=1}$ converges a.s. in ${\cal P}^B_{\b,m} G\Omega^B_{\a, m} ({\bf R}^d)$.
\end{lem}

%%%%%%%%%%%%%%%%%%%%%%%%%%%

\Proof 
In this proof the positive constant $c$ may change from line to line.
We will use Proposition \ref{pr.equi.chao} which states that, on a fixed inhomogeneous Wiener chaos, 
all $L^p$-norms ($1<p<\infty$) are equivalent.

Consider the case $(j-1)/2^k \le x \le y \le j/2^k$. 
Since $\psi(k)(t, \,\cdot\,)$ is linear on $[x,y]$, it is obvious that
$\Psi(k)^2_t (x,y) = (2^{2k}/2) (y-x)^2 \D_{j}^{k}\psi_t   \otimes \D_{j}^{k}\psi_t$.
From this and Lemma \ref{lem.est.D^2}, we see that
\begin{align}
{\mathbb E} \Bigl[  
|\Psi(k)^2_t (x,y)  - \Psi(k)^2_s (x,y) |^2
\Bigr]
&\le 
c 2^{4k} |y-x|^4 |t-s|^{\kappa /2} (2^{-k})^{1 -\kappa} 2^{-k}
\nn
\nn\\
&\le 
c |t-s|^{\kappa /2} |y-x|^{2 -\kappa}.
\label{est.psi2.1}
\end{align}

If $(j-1)/2^k \le x \le j/2^k \le y \le (j+1)/2^k$. 
From Chen's identity, we have 
\[
\Psi(k)^2_t (x,y)  
= 
\Psi(k)^2_t (x, j/2^k)  +\Psi(k)^2_t (j/2^k,y) +  \Psi(k)^1_t (x, j/2^k)  \otimes \Psi(k)^1_t (j/2^k,y).
\]
This implies 
the inequality (\ref{est.psi2.1}) in this case, too (for  different $c>0$). 
Thanks to repeated use of Chen's identity,
a similar inequality still holds even if $\Psi(k)$ is replaced with $\Psi(k+1)$.
Combining these we have 
\begin{eqnarray}
{\mathbb E} 
\Bigl[  
|\Psi(k+1)^2_t (x,y)  - \Psi(k)^2_t (x,y)  - \Psi(k+1)^2_s (x,y)  +\Psi(k)^2_s (x,y) |^2
\Bigr]
\nn\\
\le 
c |t-s|^{\kappa /2} |y-x|^{2 -\kappa}
\quad
\quad 
(\mbox{if $|y-x| \le 1/2^k$}).
\label{est.psi5}
\end{eqnarray}

%%%%%%%%%%%%%%%%%%%%%%%%%%%%%%%%%%%%%%%%%%%%%%%%%%%%\vspace{10mm}

Now we estimate the left hand side of (\ref{est.psi2.1}), when $x= I/2^k$ and $y=J /2^k$ with $I<J$.
For simplicity we will write $z^k_i = i/2^k$.
In the same way as in pp. 69--71, Lyons and Qian \cite{lq}, 
\begin{align}
\Psi(k)^2_t (z^k_I, z^k_J)
&=
\sum_{l =I+1}^{J}  \Psi(k)^2_t (z^k_{l-1}, z^k_{l})
+
\sum_{I+1 \le r <l \le J}
\Psi(k)^1_t (z^k_{r-1}, z^k_{r}) \otimes \Psi(k)^1_t (z^k_{l-1}, z^k_{l})
\nn\\
&=
\frac12  \sum_{l =I+1}^{J}   \D_{l}^{k}\psi_t \otimes  \D_{l}^{k}\psi_t 
+
\sum_{I+1 \le r <l \le J}
\D_{r}^{k}\psi_t \otimes \D_{l}^{k}\psi_t.
\nn
\end{align}
Here, we used Chen's identity and the fact that $\psi(k) (t, \,\cdot\,)$ is linear on 
each $[z^k_{l-1}, z^k_{l}]$ as a function of $x$.
Since $\Psi(k)^1_t (z^k_{I}, z^k_{J}) = \Psi^1_t (z^k_{I}, z^k_{J}) =
\sum_{l=I+1}^J  \D_{l}^{k}\psi_t$,
we have
\begin{align}
\Psi(k)^2_t (z^k_I, z^k_J)
&=
\frac12 \Psi^1_t (z^k_{I}, z^k_{J})^{\otimes 2} 
+
\frac12 
\sum_{I+1 \le r <l \le J}
\bigl( \D_{r}^{k}\psi_t \otimes \D_{l}^{k}\psi_t -   \D_{l}^{k}\psi_t \otimes \D_{r}^{k}\psi_t \bigr).
\nn
\end{align}
We can compute $\Psi(k+1)^2_t (z^{k+1}_{2I}, z^{k+1}_{2J}) 
= \Psi(k+1)^2_t (z^{k}_{I}, z^{k}_{J}) $ in the same way.
By subtracting them, we obtain
\begin{align}
\lefteqn{
\Psi(k+1)^2_t (z^k_I, z^k_J) -\Psi(k)^2_t (z^k_I, z^k_J)
}
\nn\\
&=
\frac12 
\sum_{ l= I+1 }^J
\bigl( \D_{2l -1}^{k+1}\psi_t \otimes \D_{2l}^{k+1}\psi_t 
-   \D_{2l}^{k+1}\psi_t \otimes \D_{2l -1}^{k+1}\psi_t \bigr).
\label{eq.2nd.lrrl}
\end{align}
Since $(i,i)$-components of the right hand side of (\ref{eq.2nd.lrrl}) vanish,
we have only to compute $(i,j)$-components for distinct $i, j$.
It is immediate from (\ref{eq.2nd.lrrl}) that 
\begin{align}
\lefteqn{
\Psi(k+1)^{2;i,j}_t (z^k_I, z^k_J) -\Psi(k)^{2;i,j}_t (z^k_I, z^k_J)
- 
\Psi(k+1)^{2;i,j}_s (z^k_I, z^k_J) +\Psi(k)^{2;i,j}_s (z^k_I, z^k_J)
}
\nn\\
&=
\frac12 
\sum_{ l= I+1 }^J
\bigl( \D_{2l -1}^{k+1}\psi_t^i  \D_{2l}^{k+1}\psi_t^j
-   \D_{2l-1}^{k+1}\psi_s^i \D_{2l }^{k+1}\psi_s^j  \bigr)
\nn\\
& 
\quad
-\frac12 
\sum_{ l= I+1 }^J
\bigl( \D_{2l }^{k+1}\psi_t^i  \D_{2l-1}^{k+1}\psi_t^j 
-   \D_{2l}^{k+1}\psi_s^i  \D_{2l -1}^{k+1}\psi_s^j \bigr)
\nn\\
&= 
\frac12 
\sum_{ l= I+1 }^J
( \D_{2l -1}^{k+1}\psi_t^i -   \D_{2l-1}^{k+1}\psi_s^i  )  \D_{2l}^{k+1}\psi_t^j
- \frac12 
\sum_{ l= I+1 }^J
\D_{2l -1}^{k+1}\psi_s^i
( \D_{2l }^{k+1}\psi_t^j -   \D_{2l}^{k+1}\psi_s^j )
\nn\\
& \quad
- \frac12 
\sum_{ l= I+1 }^J
( \D_{2l }^{k+1}\psi_t^i -   \D_{2l}^{k+1}\psi_s^i  )  \D_{2l -1}^{k+1}\psi_t^j
+
\frac12 
\sum_{ l= I+1 }^J
\D_{2l }^{k+1}\psi_s^i
( \D_{2l-1 }^{k+1}\psi_t^j -   \D_{2l -1}^{k+1}\psi_s^j )
\nn\\
&=:\frac12 A_1^{i,j} -\frac12 A_2^{i,j}-\frac12 A_3^{i,j} + \frac12 A_4^{i,j}.
\label{eq.2nd.defa12}
\end{align}

From the independence of $i$th and $j$th components, we have
\begin{align}
{\mathbb E} [| A_1^{i,j} |^2]  
&=
\sum_{ l,m= I+1 }^J
{\mathbb E} [ ( \D_{2l -1}^{k+1}\psi_t^i -   \D_{2l-1}^{k+1}\psi_s^i  ) 
( \D_{2m -1}^{k+1}\psi_t^i -   \D_{2m-1}^{k+1}\psi_s^i  )  ]
\cdot
{\mathbb E} [   \D_{2l}^{k+1}\psi_t^j   \D_{2m}^{k+1}\psi_t^j]
\nn\\
&=
2\sum_{I+1 \le l <m \le J}
{\mathbb E} [ ( \D_{2l -1}^{k+1}\psi_t^i -   \D_{2l-1}^{k+1}\psi_s^i  ) 
( \D_{2m -1}^{k+1}\psi_t^i -   \D_{2m-1}^{k+1}\psi_s^i  )  ]
\nn\\
&
\qquad\qquad \qquad\qquad \qquad\qquad \qquad\qquad \qquad\qquad \qquad
\times
{\mathbb E} [   \D_{2l}^{k+1}\psi_t^j   \D_{2m}^{k+1}\psi_t^j]
\nn\\
& +
\sum_{ l= I+1 }^J
{\mathbb E} [ | \D_{2l -1}^{k+1}\psi_t^i -   \D_{2l-1}^{k+1}\psi_s^i  |^2 ]
\cdot
{\mathbb E} [   |\D_{2l}^{k+1}\psi_t^j  |^2].
\label{eq.2nd.defa34}
\end{align}

%
%%%%%%

From Lemma \ref{lem.est.D^2} and $y-x= (J-I)2^{-k}$,
the second term on the right hand side of (\ref{eq.2nd.defa34}) is dominated by
\[
c (J-I) |t-s|^{\kappa /2}  2^{-(k+1)(1-\kappa)} \cdot 2^{-(k+1)} 
\le
c  |t-s|^{\kappa /2}  (y-x) 2^{-k(1-\kappa)}. 
\]

Assume that $y-x= (J-I)2^{-k} \le 1/2$.
Then, we may use Lemmas \ref{lem.cq1} and \ref{lem.cq2} to see that
the first term on the right hand side of (\ref{eq.2nd.defa34}) is dominated by
\begin{align}
\lefteqn{
c \sum_{I+1 \le l <m \le J} 
\frac{ |t-s|^{\kappa /2}2^{-2(k+1)} }{  \{ (m-l) 2^{-k} \}^{1+\kappa}}
\cdot
\frac{2^{-2(k+1)}} {(m-l) 2^{-k} }
}
\nn\\
&\le
c |t-s|^{\kappa /2} 2^{-k (2-\kappa)} \sum_{I+1 \le l <m \le J} \frac{1} {(m-l)^{2+\kappa}}
\nn\\
&\le
c |t-s|^{\kappa /2}   2^{-k (1-\kappa)} \frac{J-I}{2^k}
\le
c |t-s|^{\kappa /2}  (y-x) 2^{-k (1-\kappa)}.
\label{eq.2nd.defa56}
\end{align}
Thus, we have estimated ${\mathbb E} [| A_1^{i,j} |^2]$.
We can also estimate ${\mathbb E} [| A_p^{i,j} |^2]$ for $p=2,3,4$ in the same way.

To sum up, we have obtained the following inequality:
For $x=I 2^{-k} ,~ y=J 2^{-k}$ with
$0 \le y-x\le 1/2$, 
\begin{eqnarray}
{\mathbb E} 
\Bigl[  
|\Psi(k+1)^2_t (x,y)  - \Psi(k)^2_t (x,y)  - \Psi(k+1)^2_s (x,y)  +\Psi(k)^2_s (x,y) |^2
\Bigr]
\nn\\
\le 
c |t-s|^{\kappa /2} |y-x| \Bigl( \frac{1}{2^k} \Bigr)^{ 1-\kappa}.
\label{est.psi2.3}
\end{eqnarray}

We will check that the above estimate (\ref{est.psi2.3}) also holds when $y-x >1/2$.
In this case, divide the interval as 
$[ I 2^{-k},  J 2^{-k}]=[x,y] = [x,1/2] \cup [1/2,y]$ and use Chen's identity.
Then, we have
\begin{align}
\Psi(k)^2_t (x,y) 
&= 
\Psi(k)^2_t (x,1/2)  + \Psi(k)^2_t (1/2,y) + \Psi(k)^1_t (x,1/2) \otimes \Psi(k)^1_t (1/2,y) 
\nn\\
&=
\Psi(k)^2_t (x,1/2)  + \Psi(k)^2_t (1/2,y) + \Psi^1_t (x,1/2) \otimes \Psi^1_t (1/2,y).
\nn
\end{align}
Since a similar equality holds for $\Psi(k+1)^2$,
we see that
\begin{align}
\lefteqn{
\Psi(k+1)^2_t (x,y) - \Psi(k)^2_t (x,y) 
}
\nn\\
&= 
\Bigl\{ \Psi(k+1)^2_t (x,1/2) -  \Psi(k)^2_t (x,1/2)  \Bigr\} 
+ 
\Bigl\{ \Psi(k+1)^2_t (1/2,y)  -\Psi(k)^2_t (1/2, y) \Bigr\}.
\nn
\end{align}
We can easily see from this that (\ref{est.psi2.3}) holds (for a possibly different 
constant $c>0$) even when $y-x = (J-I)/2^k>1/2$.

%%%%%

Now, let us consider the general case $0 \le x < y \le 1$ with $y-x \ge 2^{-k}$.
Take two integers $I$ and $J$ with $I \le J$ so that 
$0 \le I 2^{-k} -x < 2^{-k}$ and $0 \le y -J 2^{-k}  < 2^{-k}$ hold.
Let us again divide the interval as $[x,y] = [x,I2^{-k}] \cup [I2^{-k},J2^{-k}] \cup [J2^{-k},y]$ 
and use Chen's identity.
Then, we have
\begin{align}
\Psi(k)^2_t (x,y) 
&= 
\Psi(k)^2_t (x,  I2^{-k}) + \Psi(k)^2_t (I2^{-k}, J2^{-k})  + \Psi(k)^2_t (J2^{-k},y) 
\nn\\
&
+ \Psi(k)^1_t (x, I2^{-k}) \otimes \Psi^1_t ( I2^{-k}, J2^{-k})
+
 \Psi(k)^1_t (x, I2^{-k}) \otimes \Psi(k)^1_t (J2^{-k},y) 
 \nn\\
  & 
  + 
   \Psi^1_t ( I2^{-k}, J2^{-k}) \otimes \Psi(k)^1_t (J2^{-k},y)
       \nn\\
        &=: B_1(k;t) +\cdots + B_6(k;t).
                \nn
\end{align}

From (\ref{est.psi5}) and (\ref{est.psi2.3}), we can easily estimate $B_i~(1 \le i \le 3)$ as follows;
\begin{eqnarray}
{\mathbb E} 
\Bigl[  
| B_i (k+1;t) -B_i (k;t) - B_i (k+1;s) +B_i (k;s) |^2
\Bigr]
\le 
c |t-s|^{\kappa /2} |y-x| \Bigl( \frac{1}{2^k} \Bigr)^{ 1-\kappa}
\label{est.psi2.4}
\end{eqnarray}
for $i=1,2,3$.

Estimates for $B_i~(4 \le i \le 6)$ are similar, but slightly more complicated.
For example, $B_4$ can be calculated as follows;
\begin{align}
\lefteqn{
B_4 (k+1;t) -B_4 (k;t) - B_4 (k+1;s) +B_4 (k;s)
}
\nn\\
&=
\Bigl\{ \Psi(k+1)^1_t (x, I2^{-k}) -   \Psi(k)^1_t (x, I2^{-k}) \Bigr\} \otimes \Psi^1_t ( I2^{-k}, J2^{-k})
\nn\\
& \quad
 -\Bigl\{ \Psi(k+1)^1_s (x, I2^{-k}) -   \Psi(k)^1_s (x, I2^{-k}) \Bigr\} 
 \otimes \Psi^1_s ( I2^{-k}, J2^{-k})
 \nn\\
  &=
   \Bigl\{ \Psi(k+1)^1_t (x, I2^{-k}) -   \Psi(k)^1_t (x, I2^{-k}) 
   \nn\\
   &\qquad \qquad
   - \Psi(k+1)^1_s (x, I2^{-k}) +   \Psi(k)^1_s (x, I2^{-k}) 
       \Bigr\} 
    \otimes \Psi^1_t ( I2^{-k}, J2^{-k}) 
\nn\\
& \quad 
- \Bigl\{ \Psi(k+1)^1_s (x, I2^{-k}) -   \Psi(k)^1_s (x, I2^{-k}) \Bigr\} 
\otimes
\Bigl\{ \Psi^1_t ( I2^{-k}, J2^{-k}) -   \Psi^1_s ( I2^{-k},J2^{-k}) \Bigr\}.
\nn
\end{align}

From (\ref{est.lmd}) and Lemma \ref{lem.est.D^2},  we see that 
\begin{align}
\lefteqn{
{\mathbb E} 
\Bigl[  
| B_4 (k+1;t) -B_4 (k;t) - B_4 (k+1;s) +B_4 (k;s) |^2
\Bigr]
}
\nn\\
&\le 
c(t-s)^{\kappa /2} ( I2^{-k}-x)^{1 -\kappa} \cdot (J-I)2^{-k}
+
c ( I2^{-k}-x) \cdot (t-s)^{\kappa /2} \{(J-I)2^{-k}\}^{1 -\kappa}
\nn\\
&\le 
c(t-s)^{\kappa /2} \Bigl( \frac{y-x}{(2^k)^{1 -\kappa} } 
+ \frac{(y-x)^{1 -\kappa} } {2^k} \Bigr).
\label{est.psi2.5}
\end{align}
If $2^{-k} \le y-x$, the right hand side is dominated by $c(t-s)^{\kappa /2} (y-x) 2^{-k(1 -\kappa)}$.
$B_5$ and $B_6$ are dominated in the same way.
Hence, the inequality (\ref{est.psi2.3}) holds for all $t, s$ and all $x, y$
with $y-x \ge 2^{-k}$.

From this and (\ref{est.psi5}), we finally obtain 
\begin{eqnarray}
{\mathbb E} 
\Bigl[  
|\Psi(k+1)^2_t (x,y)  - \Psi(k)^2_t (x,y)  - \Psi(k+1)^2_s (x,y)  +\Psi(k)^2_s (x,y) |^2
\Bigr]
\nn\\
\le 
c |t-s|^{\kappa /2} |y-x| \Bigl\{ 
 |y-x|^{1 -\kappa} \wedge  \Bigl( \frac{1}{2^{k}} \Bigr)^{1 -\kappa}
\Bigr\}
\label{est.psi2.6}
\end{eqnarray}
for all $t,s \in [0,T]$, $k \in {\bf N}$ and $0 \le x \le y \le 1$.

%%%%%%%%%%%%%%%%%%%%%%%%%%%%%%%%%%%%%%%

From (\ref{est.psi2.6})  and Proposition \ref{pr.equi.chao} in the finite dimensional case, 
we see that
\begin{align}
\lefteqn{
{\mathbb E} \Bigl[  \| \Psi (k+1)^2 - \Psi(k)^2 \|_{B; 2\b,2\a,m/2}^{m/2}  \Bigr]
}
\nn\\
&\le
c
\iint_{ {\cal S} (T)}  dsdt 
\iint_{  {\cal S} }
\frac{ 
|t-s|^{\kappa m/8} |y-x|^{m/4}
 \Bigl\{ 
 |y-x|^{(1 -\kappa)m/4} \wedge  \Bigl( \frac{1}{2^{k}} \Bigr)^{(1 -\kappa) m/4} \Bigr\}}
{|t-s|^{1 +\b m} |y-x|^{1 +\alpha m} 
}
dxdy.
\nn
\end{align}
If $\kappa \in (0,1)$ is chosen so that $4\b <\kappa /2 <1-2\a$,
then the integral on the right hand side converges.
Moreover, we can see from (\ref{eas.int.form}) that the inequality 
 (\ref{ineq.Bcon2}) holds with $\eta = 2^{-  \{ 1 -2 \a -(\kappa /2) \} } \in (0,1)$.
From the inequality in the lemma, almost sure and $L^{m/2}$ convergences 
of $\{\Psi (k)^2 \}$ are immediate.
Using Proposition \ref{pr.equi.chao},
we finish the proof of the lemma.
\QED

As a corollary, we obtain convergence of 
$\{\Psi (k) \}$
%$= (\psi(\cdot,0), \Psi (k)^1, \Psi (k)^2)\}$ 
%
in 
${\cal P}_{\infty} G\hat\Omega^H_{\a} ({\bf R}^d)$
by using Besov-H\"older embedding.
%
%\end{sloppypar}

\begin{cor}
Assume $\a \in (1/3, 1/2)$.
Then, 
$\{\Psi (k) \}_{k=1}^{\infty}$ converges to $\Psi$ a.s. as $k \to \infty$ 
in ${\cal P}_{\infty} G\hat\Omega^H_{\a} ({\bf R}^d)$.
Moreover,
\[
\sup_{0 \le t \le T} \| \Psi (k)^1_t -\Psi^1_t \|_{H;\a}  
+
\sup_{0 \le t \le T} \| \Psi (k)^2_t -\Psi^2_t \|_{H;2\a}  
\]
converges to $0$ in $L^p$ as $k \to \infty$ for all $1<p<\infty$.
\end{cor}

\Proof
Take $\hat\a \in (\a,1/2)$. 
Choose $\b>0$ and  $m>2$ so that $\b >1/m$, $\hat\a -1/m >\a$ and $4\b <1-2 \hat\a$.
Then, the corollary is immediate from Lemma \ref{lem.Bcon1}, 
Lemma \ref{lem.Bcon2}, and  Propoition \ref{pr.PGL.inj} (in particular, Eq. (\ref{est_FV.A++})).
\QED

%\newpage
%%%%%%%%%%%%%%%%%%%%%%%%%%%%%%%%%%%%%%%%%%%%%%%%%%%%%%%%%%%%%%%%%%%%%%%%%%%%%%%%%%%
%%%%%%%%%%%%%%%%%%%%%%%%%%%%%%%%%%%%%%%%%%%%%%%%%%%%%%%%%%%%%%%%%%%%%%%%%%%%%%%%%
%%      Section    LDP
%%%%%%%%%%%%%%%%%%%%%%%%%%%%%%%%%%%%%%%%%%%%%%%%%%%%%%%%%%%%%%%%%%%%%%%%%%%%%%%%%%%
%%%%%%%%%%%%%%%%%%%%%%%%%%%%%%%%%%%%%%%%%%%%%%%%%%%%%%%%%%%%%%%%%%%%%%%%%%%%%%%%%%%

\section{Large deviation principle}\label{sec.ldp}
%\noindent
%$\spadesuit$~Cameron-Martin space.
%
%

In this section we will state and prove our main theorem (Theorem \ref{pr.ldp.st}),
using the method developed by Friz and Victoir in \cite{fv07}.
First, we will study the regularity of a generic element of 
Cameron-Martin space.
Next, we will prove a large deviation principle for a weaker topology.
Finally, we will strengthen the topology 
by using the inverse contraction principle and  exponential tightness.

We denote by $\mu$ the law of the two parameter process 
$\psi= (\psi(t,x))_{0 \le t \le T, 0 \le x \le 1}$.
For $\a \in (1/3, 1/2)$, $\b >0$, and $m \ge 1$ such that $\a -1/m >1/3$, $\b >1/m$, and
$4\b <1-2\a$,
$\mu$ is a Gaussian measure on 
$C_0^{B;\b,m} ( [0,T], {\cal X} )$ with ${\cal X} = C^{B;\a,m} ( [0,1], {\bf R}^d)$.
We denote by ${\cal H}$
Cameron-Martin space of $\mu$.

First we investigate the regularity of $x \mapsto h(t,x)$
for $h \in {\cal H}$ to make sure that $h$ admits a spatial lift ${\cal L}_2 (h)$.
\begin{lem}\label{lem.regH}
Let $h \in {\cal H}$ and $t \in [0,T]$.
\\
\noindent
{\rm (i)}~
$h(t, \cdot)$ is $1/2$-H\"older continuous and its
$1/2$-H\"older norm is dominated by $c\|h\|_{\cal H}$ for some positive 
constant $c=c_T$ independent of $h$ and $t$.
\\
\noindent
{\rm (ii)}~
$h(t, \cdot)$ is of finite $q$-variation for any $q \in (4/3, 2)$
and its $q$-variation norm is dominated by $c\|h\|_{\cal H}$ for some positive 
constant $c=c_{T,q}$ independent of $h$ and $t$.
\end{lem}

\Proof
We may assume $d=1$.
Let $Z$ be the measurable linear functional (i.e., the element in the first order Wiener chaos) 
associated with $h$.
$Z$ is a mean-zero Gaussian random variable with variance $\|h\|_{{\cal H}}^2$. 
Then, from the general theory of abstract Wiener spaces and the fact that 
the evaluation map at $(t,x)$ is a continuous linear functional, 
we see that ${\mathbb E}[Z \psi(t,x)] =h(t,x)$.
Here, ${\mathbb E}$ stands for the expectation with respect to $\mu$.

By Corollary \ref{cor.kolm} and Schwarz' inequality,  we can easily see that
\[
|h(t,y) - h(t,x)| \le {\mathbb E}[Z^2]^{1/2} {\mathbb E}[|\psi(t,y) - \psi(t,x)|^2] 
\le 
c \|h\|_{{\cal H}} |t-s|^{1/2}.
\]
Thus we have proved the first assertion.

Now we show the second assertion, following p. 783, Friz and Victior \cite{fv07}.
Recall that Wick's formula implies that 
\[
{\mathbb E}[U^2 V^2] 
=
2{\mathbb E}[U V]^2  + {\mathbb E}[U^2] {\mathbb E}[V^2] 
\qquad
\mbox{or}
\qquad
{\rm Cov} (U^2, V^2) = 2{\mathbb E}[U V]^2 
\] 
for two mean-zero Gaussian random variables $U, V$.
Then 
\begin{eqnarray*}
Q_k :&=& 
\sum_{i=1}^{2^k} |h(t, \frac{i}{2^{k}}) - h(t, \frac{i-1}{2^{k}} ) |^2
=
\sum_{i=1}^{2^k} 
{\mathbb E}[Z  \Psi^1_t (\frac{i-1}{2^{k}}, \frac{i}{2^{k}})   
]^2 
\nn\\
&=&
\Bigl\{ \sum_{i=1}^{2^{k-1}} + \sum_{i=2^{k-1} +1}^{2^k}  \Bigr\}
{\mathbb E}[Z  \Psi^1_t (\frac{i-1}{2^{k}}, \frac{i}{2^{k}})   
]^2 
=: Q_k^{(1)} +Q_k^{(2)}.
\end{eqnarray*}
The first term $Q_k^{(1)}$ on the right hand side is dominated as follows:
\begin{align}
Q_k^{(1)} 
&= 
\frac12 \sum_{i=1}^{2^{k-1}} {\mathbb E} \bigl[ Z^2 
 \bigl\{ \Psi^1_t (\frac{i-1}{2^{k}}, \frac{i}{2^{k}})^2 -  
 {\mathbb E} [ \Psi^1_t (\frac{i-1}{2^{k}}, \frac{i}{2^{k}})^2]
   \bigr\}    
\bigr] 
\nn\\
&\le
\frac12  {\mathbb E}[Z^4]^{1/2}
\Bigl\{  
\sum_{i, l=1}^{2^{k-1}} {\rm Cov} \Bigl( \Psi^1_t (\frac{i-1}{2^{k}}, \frac{i}{2^{k}})^2,
\Psi^1_t (\frac{l-1}{2^{k}}, \frac{l}{2^{k}})^2\Bigr)
\Bigr\}^{1/2}
\nn\\
&=
c \|h\|_{{\cal H}}^2 
\Bigl\{  
\sum_{i, l=1}^{2^{k-1}} 
{\mathbb E} \bigl[
\Psi^1_t (\frac{i-1}{2^{k}}, \frac{i}{2^{k}})  \Psi^1_t (\frac{l-1}{2^{k}}, \frac{l}{2^{k}})
\bigr]^2
\Bigr\}^{1/2}.
\nn
\end{align}

From Corollary \ref{cor.kolm} and Lemma \ref{lem.cq1}, we have
\begin{align}
\lefteqn{
\sum_{i, l=1}^{2^{k-1}} 
{\mathbb E} \bigl[
\Psi^1_t (\frac{i-1}{2^{k}}, \frac{i}{2^{k}})  \Psi^1_t (\frac{l-1}{2^{k}}, \frac{l}{2^{k}})
\bigr]^2
}\nn\\
&
=
\Bigl\{
\sum_{1 \le  i,l \le 2^{k-1} :|i-l| \le 1} + \sum_{1 \le  i,l \le 2^{k-1} :|i-l| \ge 2}
\Bigr\}
{\mathbb E} \bigl[
\Psi^1_t (\frac{i-1}{2^{k}}, \frac{i}{2^{k}})  \Psi^1_t (\frac{l-1}{2^{k}}, \frac{l}{2^{k}})
\bigr]^2
\nn\\
&\le
c 2^{k-1} (2^{-k})^2 + c \sum_{1 \le  i,l \le 2^{n-1} :|i-l| \ge 2}
\Bigl| \frac{(2^{-k})^2}{ i2^{-k} -l2^{-k}} \Bigr|^2
\le c 2^{-k}.
\nn
\end{align}
Here, we have used $\sum_{l=1}^{\infty} l^{-2} <\infty$.
The second term $Q_k^{(2)}$ can be dominated in the same way.
Consequently, $Q_k \le c  \|h\|_{{\cal H}}^2  2^{-k/2}$.

We use the following inequality (see Proposition 4.1.1 in p. 62, \cite{lq}):
for $\gamma > q-1$, there is a positive constant $c=c_{q,\gamma}$ such that
$$
\| h(t, \,\cdot\,) \|_{var ;q}^q 
\le 
c_{q,\gamma} \sum_{k=1}^{\infty} k^{\gamma} 
\sum_{i=1}^{2^k}
 |h(t, \frac{i}{2^{k}}) - h(t, \frac{i-1}{2^{k}} ) |^q.
$$
The right hand side is dominated by 
$$   
   c \sum_{k=1}^{\infty} k^{\gamma} Q_k^{q/2} 2^{k (1- q/2)}
     \le 
    c  \|h\|_{{\cal H}}^q \sum_{k=1}^{\infty}  k^{\gamma} 2^{k \{ 1 - (3q)/4 \}},
$$
which is convergent if $q > 4/3$.
This proves the second assertion of Lemma \ref{lem.regH}.
\QED

By Lemma \ref{lem.regH}, the spatial lift  $H={\cal L}_2 (h)$ exists for all $h$.
A natural question is whether and in which topology   $H(k)={\cal L}_2  (h(k))$
converges to ${\cal L}_2 (h)$.
Actually, it seems rather difficult to prove the convergence with respect to H\"older or Besov topology.
However, it is much easier to prove it with respect to a weaker topology,
namely,
the uniform topology.
\begin{lem}\label{lem.uconvH}
For any $L>0$ and $i=1,2$,
\begin{eqnarray*}
\sup\bigl\{ 
| H(k)^i_t(x,y) - H^i_t(x,y)|
\,\mid  \,
(x, y) \in {\cal S}, t \in [0, T],
\,
\mbox{ $h \in {\cal H}$ with $\|h\|_{{\cal H}} \le L$}
\bigr\}
\end{eqnarray*}
converges to zero as $k \to \infty$.
\end{lem}

\Proof
From Lemma \ref{lem.regH} {\rm (i)}, 
\[
|h(t, x) - h(k)(t,x)| \le 2 \|h(t, \,\cdot\,)\|_{H; 1/2} (2^{-k})^{1/2} \le cL 2^{-k/2}.
\]
The constant $c$ on the right hand side is independent of $t,x,k, L$.
Taking a difference we can easily prove the case $i=1$.

Next we consider the case $i=2$.
Let $4/3 < q <r <2$.
From  Lemma 1.12 in p. 8, Lyons, Caruana, and Levy \cite{lcl}
and Lemma \ref{lem.regH} {\rm (ii)} above, we see that
\[
\|h(k)(t, \,\cdot\,)\|_{var;q} \le \|h(t, \,\cdot\,)\|_{var;q} \le c L 
\]
and
\begin{eqnarray*}
\lefteqn
{
\|h(k)(t, \,\cdot\,) - h(t, \,\cdot\,) \|_{var;r}^r
}
\nn\\
&\le
c
\|h(k)(t, \,\cdot\,) - h(t, \,\cdot\,) \|_{\infty}^{r-q} 
\{
\|h(k)(t, \,\cdot\,) \|_{var;q}^q + \|h(t, \,\cdot\,) \|_{var;q}^q 
\}
\le
c L^r 2^{- k (r-q)/2}.
\end{eqnarray*}
Therefore, the left hand side converges to zero uniformly in $t$ and $h$ with $\|h\|_{{\cal H}} \le L$.
By the theory of Young integral, we have
\begin{align}
\lefteqn{
|H(k)^2_t(x,y) - H^2_t(x,y)|
}
\nn\\
&\le
\bigl|
\int_x^y \{ H(k)^1_t(x,z) - H^1_t(x,z) \} \otimes d_z h(k)(t,z)
\bigr|
\nn\\
&\qquad \qquad+
\bigl|
\int_x^y  H^1_t(x,z)  \otimes d_z\{ h(k)(t,z) -h(t,z) \}
\bigr|
\nn\\
&\le 
c \|h(k)(t, \,\cdot\,) - h(t, \,\cdot\,) \|_{var;r}
\{ 
\|h(k)(t, \,\cdot\,) \|_{var;r} +\| h(t, \,\cdot\,) \|_{var;r}
\}
\le
cL^2 2^{- k (1-q/r)/2}. 
\nn
\end{align}
The constant $c$ on the right hand side is independent of $t,x,y, k, L$.
Letting $k \to \infty$, we prove the case $i=2$.
\QED

%%%%%%%%%%%%%%%%%%%%%%%%%%%%%%%%%%%%%%%%%%%%%%%%%%%%%%%%%%%%%%%%%%%%%%%%%%

We will show below that $\{ \ve \Psi (k)\}$ is
an exponentially good approximation of $ \ve \Psi$ with respect to the uniform topology.
\begin{lem}\label{lem.expgood}
There exists a distance ${\rm dist}_{\infty}$  on ${\cal P}_{\infty} G \hat\Omega_{\infty} ({\bf R}^d)$
which induces the same uniform topology on this set and satisfies the following property: 
For any $\delta >0$, it holds that
\[
\lim_{k \to \infty}  \limsup_{\ve \searrow 0 }
\ve^2 
\log {\mathbb P} \bigl(  {\rm dist}_{\infty} ( \ve \Psi, \ve \Psi (k)  ) > \delta \bigr)   = -\infty.
\]
\end{lem}

\Proof
A geometric rough path is a continuous path that takes values in 
the free nilpotent group $G^2 =G^2 ({\bf R}^d)$ of step 2.
We take this viewpoint, following p. 776, \cite{fv07} and Chapters 7--9, \cite{fvbk}. 
$G^2$
 is a subset of $\{ (1, g_1, g_2) ~|~ g_1 \in {\bf R}^d,
 g_2\in ({\bf R}^d)^{\otimes 2}   \}\subset T^{2}({\bf R}^d)$.
The multiplication of $G^2$ is the tensor product $\otimes$ in $T^{2}({\bf R}^d)$,
 which will be suppressed.
On $G^2$
there exists a homogeneous and symmetric norm $\| \,\cdot\, \|$.
It satisfies that $ \| g \| = \| g^{-1} \|$
and $\| \lm g \| = |\lm| \| g \|$ for the dilation by $\lm \in {\bf R}$.
Moreover, there exists a constant $c \ge 1$ such that
\[
c^{-1}\| g \| 
\le
|g_1|_{{\bf R}^d}  + |g_2|_{({\bf R}^d)^{\otimes 2}}^{1/2}
\le
c\|  g \| 
\]
for all $g=(1,g_1, g_2) \in G^2$.
A distance ${\rm d}_{G^2}$ is defined by 
${\rm d}_{G^2} (g, \hat{g}) := \| g^{-1} \hat{g} \| =  \| \hat{g}^{-1} g\|$.

It is known that, for $A \in G\Omega_{\infty}({\bf R}^d)$, 
$[0,1] \ni x \mapsto (1,A^1_{0,x}, A^2_{0,x})$
is a continuous path in $G^2$ which starts at the unit element.
Conversely, for a  continuous path $a$ in $G^2$ which starts at the unit element,
$A_{x,y} = a_x^{-1} a_{y}$ defines an element in $G\Omega_{\infty}({\bf R}^d)$.
Hence, these spaces sets $G\Omega_{\infty}({\bf R}^d)$ 
and $C_{{\rm unit}}([0,1], G^2)$ can be identified as sets.
Topology of uniform convergence is induced by the following distance;
$\hat{\rm d}(a,b) = \sup_{0 \le x \le 1} {\rm d}_{G^2}(a_x, b_x)$.
It is easy to see that for some constant $c \ge 1$ 
\[
c^{-1} \hat{\rm d}(a,b)
\le 
\sup_{0 \le x \le 1}
\Bigl\{
|B^1_{0,x} - A^1_{0,x} |+ |(B^2_{0,x} - A^2_{0,x}) - A^1_{0,x} \otimes (B^1_{0,x} - A^1_{0,x}) |^{1/2}
\Bigr\}
\le 
c\hat{\rm d}(a,b)
\]
Note that $\lim_{k \to \infty} A(k) = A$ in $G\Omega_{\infty}({\bf R}^d)$ 
(that is, $\lim_{k \to \infty} \sup_{(x,y) \in {\cal S}} |A(k)^i_{x,y} - A^i_{x,y}|$
for $i=1,2$)
is equivalent to $\lim_{k \to \infty}\hat{\rm d}(a(k), a) =0$.

Let us consider ${\cal P}_{\infty} G\Omega_{\infty}({\bf R}^d)$,
which can be identified with $C([0,T], G\Omega_{\infty}({\bf R}^d))$ as sets.
Uniform convergence on this set is induced by the following distance;
$$
{\rm dist}_{\infty} ( a,b) = \sup_{0 \le t \le T} \hat{\rm d}(a_t,b_t)
=
\sup_{0 \le t \le T} \sup_{0 \le x \le 1} {\rm d}_{G^2}(a_{t,x}, b_{t,x}),
$$
where $a,b \in C([0,T], G\Omega_{\infty}({\bf R}^d)) \cong C([0,T],  C_{{\rm unit}}([0,1], G^2))$.
For some constant $c \ge 1$,  
\begin{align}
%\lefteqn{
c^{-1}   {\rm dist}_{\infty} (a,b)
%}\nn\\
&\le 
\sup_{0 \le t \le T} \sup_{0 \le x \le 1}
\Bigl\{
|B^1_t(0,x) - A^1_t(0,x) |
\nn\\
& \quad
+ \bigl|(B^2_t(0,x) - A^2_t(0,x)) - A^1_t(0,x) \otimes (B^1_t(0,x) - A^1_t(0,x)) \bigr|^{1/2}
\Bigr\}
\nn\\
&\le 
c {\rm dist}_{\infty} (a,b).
\label{est.dinfty}
\end{align}
Again
 $\lim_{k \to \infty} \sup_{t \in [0,T], (x,y) \in {\cal S}} |A(k)^i_t(x,y) - A^i_t(x,y)|$
for $i=1,2$
is equivalent to convergence with respect to ${\rm dist}_{\infty}$.
This distance naturally extends to a distance on ${\cal P}_{\infty} G \hat\Omega_{\infty} ({\bf R}^d)$
with a trivial modification
and is again denoted by ${\rm dist}_{\infty}$.
The reason why we introduce this distance is 
that ${\rm dist}_{\infty} (\ve a, \ve b) =\ve {\rm dist}_{\infty} (a,b)$
for the dilation by $\ve \ge 0$.

Take any $\a, \b, m$ which satisfy the assumptions of 
Lemmas \ref{lem.Bcon1} and \ref{lem.Bcon2}.
By Proposition \ref{pr.PGL.inj}, $\Psi(k)$ also converges to $\Psi$
in ${\cal P}_{\infty} G\Omega_{\infty}({\bf R}^d)$ a.s.
Moreover, we see from (\ref{est_FV.A++}) and Lemmas \ref{lem.Bcon1}, \ref{lem.Bcon2}  that
 $ 
 \sup_{t \in [0,T], (x,y) \in {\cal S}} |\Psi (k)^i_t(x,y) - \Psi^i_t(x,y)|
  $
  converges to $0$ in any $L^p$-norm ($1 <p <\infty$).

%\begin{sloppypar}
Let $s_k$ be the maximum of 
{\rm (i)}~$L^4$-norm of $\sup_{t, (x,y)} |\Psi (k)^1_t(x,y) - \Psi^1_t(x,y)|$,
{\rm (ii)}~$L^2$-norm of $\sup_{t, (x,y)} |\Psi (k)^2_t(x,y) - \Psi^2_t(x,y)|$,
and  {\rm (iii)}~$L^2$-norm of 
$\sup_{t, (x,y)}  | \Psi^1_t(x,y) \otimes (\Psi^1_t(x,y) - \Psi(k)^1_t(x,y)) |$.
Then, $\lim_{k \to \infty} s_k =0$.
%\end{sloppypar}
%
%
%
By Proposition \ref{pr.equi.chao}, we have
\begin{align}
 \Bigl\| {\rm dist}_{\infty} (\Psi, \Psi(k) )   \Bigr\|_{L^{4q}}
&\le 
c\Bigl\| 
\sup_{t, (x,y)} |\Psi (k)^1_t(x,y) - \Psi^1_t(x,y)|
 \Bigr\|_{L^{4q}}
 \nn\\
 &\quad + c
   \Bigl\| 
\sup_{t, (x,y)} |\Psi (k)^2_t(x,y) - \Psi^2_t(x,y)|
 \Bigr\|_{L^{2q}}^{1/2}
  \nn\\
 &\quad +c
   \Bigl\| 
   \sup_{t, (x,y)}  | \Psi^1_t(x,y) \otimes (\Psi^1_t(x,y) - \Psi(k)^1_t(x,y)) |
      \Bigr\|_{L^{2q}}^{1/2}
       \le 
           c' s_k \sqrt{4q}.
     \nn
    \end{align}
Here, $c, c'$ are positive constants independent of $q \ge 1$ and  $k= 1,2,\ldots$.

Therefore, for sufficiently large $q$, we have
\begin{align}
{\mathbb P} \bigl(  {\rm dist}_{\infty} ( \ve \Psi, \ve \Psi (k)  ) > \delta \bigr) 
&=
{\mathbb P} \bigl(  {\rm dist}_{\infty} (  \Psi,  \Psi (k)  ) > \delta /\ve \bigr) 
\nn\\
&\le 
(\delta /\ve)^{-q} c'^q  s_k^q \sqrt{q}^q
\nn\\
&= \exp \Bigl[  
q \log \bigl( \frac{\ve c'}{\delta}  s_k \sqrt{q}  \bigr)
\Bigr].
\nn
\end{align}
Choosing $q =1/\ve^2$, we have
  \[
\ve^2 \log  {\mathbb P} \bigl(  {\rm dist}_{\infty} ( \ve \Psi, \ve \Psi (k)  ) > \delta \bigr) 
  \le 
  \log \bigl( \frac{c'}{\delta}  s_k \bigr).
  \]
First taking $\limsup_{\ve \searrow 0}$ and then letting $k \to \infty$,
we complete the proof.  \QED

%%%%%%%%%%%%%%%%%%%%%%%%%%%%%%%

%
Using Lemmas \ref{lem.uconvH} and \ref{lem.expgood}, we can  prove 
a Schilder-type large deviation principle for $\nu_{\ve}$, which is 
the laws of $\ve \Psi$,
with respect to a weaker topology.
Let us define a rate function $I$ by
$I (A) = \|h\|_{{\cal H}}^2 /2$ if $A = {\cal L}_2(h)$
for some $h \in {\cal H}$
and 
$I (A) =\infty$ if 
$A \in {\cal P}_{\infty} G \hat\Omega_{\infty} ({\bf R}^d) \setminus {\cal L}_2 ({\cal H})$.
Here, ${\cal L}$ denotes the natural lift with respect to the spatial parameter.

%%%%%%%%%%%

\begin{prop}\label{pr.ldp.weak}
Let the notations be as above. 
Then, 
$(\nu_{\ve})_{0<\ve \le 1}$ satisfies a large deviation principle 
as $\ve \searrow 0$ in ${\cal P}_{\infty} G \hat\Omega_{\infty} ({\bf R}^d)$
with a good rate function $I$.
\end{prop}

\Proof
The laws of $\ve \psi$ induces scaled Gaussian measures on $C_0 ([0,T], C([0,1],{\bf R}^d))$.
By the general theory of abstract Wiener spaces, 
they satisfy Schilder's large deviation with a good rate function $J$.
Here, $J(h) = \|h\|_{{\cal H}}^2 /2$ if $h \in {\cal H}$
and $J(h)=\infty$ if $h \notin {\cal H}$.
The map $a \to A(k) ={\cal L}_2(a(k))$
is clearly continuous
for each fixed $k=1,2,\ldots$.
Therefore, the laws of $\ve \Psi(k)$ satisfies a large deviation with a good rate function
since it is a continuous image of $\ve \psi$ (see Theorem 4.2.1, \cite{dz}).

From Lemma \ref{lem.uconvH} and (\ref{est.dinfty}),  it is immediate that 
$\sup_{ \|h\|_{{\cal H}} \le L } {\rm dist}_{\infty} (H(k), H) \to 0$
as $k \to 0$ for each $L>0$.
Combining this with Lemma \ref{lem.expgood}, 
we can use exponentially good approximation (Theorem 4.2.23, \cite{dz})
to finish the proof.
\QED

%%%%%%%%%%%%%%%%%%%%%%%%%%%%%%%%%%%%%%%%%%%%%%%%%%%%%%%%%%%%%%%%%%%%%%%%%%

%\vspace{10mm}
%\noindent
%$\spadesuit$~ Exponential tightness.

We prove that the laws of $\ve \Psi$ are exponentially tight on 
the path space over the geometric rough path space.
To strengthen the topology, exponential tightness is the key.

\begin{lem}\label{lem.e^approx}
$(\nu_{\ve})_{0<\ve \le 1}$ are exponentially tight on 
${\cal P}_{\infty} G \hat\Omega^H_{\a} ({\bf R}^d)$ for any $\a \in (1/3,1/2)$.
\end{lem}

\Proof
Take $\a' \in (\a, 1/2)$ and $\b >0$ so that $4\b < 1-2\a'$.
Take $m > 2$ sufficiently  large so that $\a' -1/m >\a$ and $\b >1/m$.
Then,  Lemmas \ref{lem.Bcon1} and \ref{lem.Bcon2} hold.

Set $Z = \| \psi ( \,\cdot\,,0) \|_{B;\b,m } +\| \Psi^1 \|_{B;\b,\a',m }$ 
or $Z = \| \Psi^2 \|_{B; 2\b, 2\a',m/2 }^{1/2}$.
By Proposition \ref{pr.equi.chao},
$\|Z\|_{L^q} \le c \sqrt{q}$ for some $c>0$.
This implies that ${\mathbb E} [ e^{ \eta Z^2 } ] <\infty$  if $\eta \in (0,   (2c^2e)^{-1})$.
Indeed, 
\begin{align}
{\mathbb E} [ e^{ \eta Z^2 } ]  
&=
\sum_{k=0}^{\infty }  \frac{ \eta^k {\mathbb E}[ Z^{2k}] } { k! }
\le 
\sum_{k=0}^{\infty }  \frac{ \eta^k  c^{2k} (2k)^{k} } { k! },
\nn
\end{align}
which converges as a power series of $\eta$ if $|\eta| \le 1/(2c^2e)$.
By Chebyshev's inequality, this square exponential integrability implies that
for some constant $c>0$  it holds that
${\mathbb P} (Z >r) \le c e^{-\eta r^2}$ for all $r>0$.

For $r >0$,  we set
\begin{align}
B_r :&= \{   (a, A^2) =  (a_0, A^1,A^2) \in  {\cal P}_{\infty} G \hat\Omega^H_{\a} ({\bf R}^d) 
~ \mid~
\nn\\
&
\qquad\qquad
\| a ( \,\cdot\,,0) \|_{B;\b,m }+
\| A^1 \|_{B;\b,\a',m } \le r, \, \| A^2 \|_{B; 2\b, 2\a',m/2 }^{1/2} \le r
\}.
\nn
\end{align}
From Proposition \ref{pr.PGL.inj}, its closure $\bar{B_r}$ is compact in
 ${\cal P}_{\infty} G \hat\Omega^H_{\a} ({\bf R}^d) $.
The complement set $(\bar{B_r})^c$ is contained in 
%
%
%\begin{align}
$
\{ 
%  (a, A^2) \in  {\cal P}_{\infty} G \hat\Omega^H_{\a} ({\bf R}^d) 
%\,\mid\,
%\nn\\
\| a( \,\cdot\,,0) \|_{B;\b,m } +\| A^1\|_{B;\b,\a',m } > r  
\, \mbox{ or }  \,  \| A^2 \|_{B;2\b, 2\a',m/2 }^{1/2} > r\}.
$
%\nn
%\end{align}
%
%
%
Noting that $\|(\ve A)^2 \|_{B;2\b, 2\a',m/2 }^{1/2} = \ve \| A^2 \|_{B;2\b, 2\a',m/2 }^{1/2}$,
we can easily see that
\[
\limsup_{\ve \searrow 0} \ve^2 \log \nu_{\ve} ( (\bar{B_r})^c )   \le -\eta r^2.
\]
Letting $r \to \infty$, we prove  exponential tightness.
\QED

%%%%%%%%%%%%%%%%%%%%%%%%%%%%%%%%%%%%%%%%%%%%%%%%%%%%%%%%%%%%%%%%%%%%%%%%%%

%\vspace{10mm}
%\noindent
%$\spadesuit$~ 
%LDP with respect to  stronger topology.

Now we state and prove our main theorem.
Thanks to exponential tightness in Lemma \ref{lem.e^approx},
we can prove a large deviation principle of Schilder type
with respect to the desired topology using the inverse contraction principle.
It is a little bit interesting that we can show
${\cal L}_2 ({\cal H}) \subset {\cal P}_{\infty} G \hat\Omega^H_{\a} ({\bf R}^d)$ 
in a rather indirect way like this.

\begin{thm}\label{pr.ldp.st}
Let $\a \in (1/3,1/2)$
and let $\nu_{\ve}$ be the law of $\ve \Psi$. 
Then,
$(\nu_{\ve})_{0< \ve \le 1}$ satisfies a large deviation principle 
as $\ve \searrow 0$ in ${\cal P}_{\infty} G \hat\Omega_{\a}^H ({\bf R}^d)$
with a good rate function $I$.
In particular,  the effective domain of $I$, i.e., ${\cal L}_2 ({\cal H})$,
is  a subset of ${\cal P}_{\infty} G \hat\Omega^H_{\a} ({\bf R}^d)$.
\end{thm}

\Proof
By Lemma \ref{lem.e^approx}, we can use the inverse contraction principle.
(See Theorem 4.2.4 and Remarks in pp.127--129, \cite{dz}.)
The inclusion ${\cal L}_2 ({\cal H}) \subset {\cal P}_{\infty} G \hat\Omega^H_{\a} ({\bf R}^d)$
is a part of the inverse contraction principle, 
but we will prove it now for the reader's convenience.
Suppose ${\cal L}_2 (h) \notin {\cal P}_{\infty} G \hat\Omega^H_{\a} ({\bf R}^d)$
for some $h \in {\cal H}$.
By exponential tightness, there exists a compact subset 
$K \subset {\cal P}_{\infty} G \hat\Omega^H_{\a} ({\bf R}^d)$
such that $\limsup_{\ve \searrow 0} \ve^2 \log \nu_{\ve}(K^c) < - (1+ \|h\|^2_{{\cal H}} /2 )$.
Since the injection is continuous (Proposition \ref{pr.PGL.inj}),
$K$ is also compact in ${\cal P}_{\infty} G \hat\Omega_{\infty} ({\bf R}^d)$ and,
in particular, $K^c$ is open.
By Proposition \ref{pr.ldp.weak}, 
$\liminf_{\ve \searrow 0} \ve^2 \log \nu_{\ve}(K^c) \ge - \|h\|^2_{{\cal H}} /2$,
since ${\cal L}_2 (h) \in K^c$. 
This is a contradiction.
\QED

%%%%%%%%%%%%%%%%%%%%%%%%%%%%%%%%%%%%%%%%%%%%%%%%%%%%%%%%%%%%%%%%%%%%%%%%%

%\vspace{10mm}
%\noindent
%$\spadesuit$~ 

Before we end this section,
we give a comment on Freidlin-Wentzell type large deviations.
By the well-known contraction principle (Theorem 4.2.1, \cite{dz}),
the laws of $F(\ve \Psi)$ also satisfies a large deviation principle
for any continuous map $F$ from 
${\cal P}_{\infty} G \hat\Omega_{\a}^H ({\bf R}^d)$ to any  Hausdorff space.

%In the additive case of M. Hairer's rough  stochastic PDE theory, 
%the solution map is a continuous map of $\Psi$ with respect to the topology of 
%${\cal P}_{\infty} G \hat\Omega_{\a}^H ({\bf R}^d)$.
%Hence, the laws of the solutions of the rough stochastic PDE 
%automatically satisfy a large deviation principle of Freidlin-Wentzell type.
%Below we will give a short remark on it below.
%
%

Hairer and Weber \cite{hw} considered the following rough stochastic PDE;
\begin{equation}
d_t u = \bigl[
\triangle u + f(u) + g(u) \partial_x u 
\bigr] dt + \theta (u)  dW_t,
\qquad
\mbox{with $u(0,x)=u_0(x)$.}
\label{eq.spde1}
\end{equation}
Here, 
{\rm (i)}~ the solution $u =u(t,x)$ is a function  from $[0,T] \times S^1$ to ${\bf R}^d$;
{\rm (ii)}~ the coefficients  
$f: {\bf R}^d \to {\bf R}^d$,
$g: {\bf R}^d \to ({\bf R}^d)^{\otimes 2}$, 
$\theta: {\bf R}^d \to ({\bf R}^d)^{\otimes 2}$, are sufficiently regular functions; 
{\rm (iii)}~$u_0(x)$ is a given initial condition at $t=0$;
{\rm (iv)}~$(W_t)_{0 \le t \le T}$ is an $L^2 (S^1, {\bf R}^d)$-cylindrical Brownian motion.
Note that $u =\psi$
when $f \equiv 0\equiv g$ and $\theta \equiv {\rm Id}$.

In the additive noise case, i.e., $\theta \equiv {\rm Id}$, 
it is proved in \cite{hw} that under mild assumptions on 
the coefficients and the initial condition, stochastic PDE (\ref{eq.spde1}) 
is well-defined in rough path sense and has a unique global solution.
Moreover, $u$ is a continuous image of $\Psi$ with respect to 
${\cal P}_{\infty} G \hat\Omega^H_{\a} ({\bf R}^d)$-topology for some $\a \in (1/3, 1/2)$.
(For precise information on the sufficient conditions for $f, g, u_0$ 
and on the topology of the space in which the solution $u$ lives, see \cite{hw, hmw}, etc.)

As a result, we see from our main theorem and the contraction principle that
a large deviation principle of Freidlin-Wentzell type holds in this case 
for the laws of the solution $u^{\ve}=u^{\ve} (t,x)$ of the following rough SPDE;
\begin{equation}
d_t u^{\ve} = \bigl[
\triangle u^{\ve} + f(u^{\ve}) + g(u^{\ve}) \partial_x u^{\ve} 
\bigr] dt +   \ve dW_t,
\qquad
\mbox{with $u^{\ve}(0,x)=u_0(x)$.}
%\label{eq.spde2}
\nn
\end{equation}
Freidlin-Wentzell type large deviations for various stochastic PDEs were extensively studied.
%For example, see \cite{cr, cm, ch, dm, fw, gs, kx, pe, so, za, zh} among many others.
In none of them was rough path theory used, however.
(Unfortunately, in the multiplicative case, i.e., the case $\theta$ is not a constant, 
the continuity of the map is not known.
So we cannot use the contraction priciple to 
prove Freidlin-Wentzell type large deviations at present.)

\begin{remark}
We studied (the lift of) the solution of stochastic heat equation (\ref{eq.psi.spde}).
In Hairer's theory, however, two-parameter Gaussian processes 
differ from paper to paper.
See \cite{fggr, hai1, hai2}.
Therefore, it may be interesting to study whether similar large deviations hold 
for the lifts of these processes, too.
\end{remark}

%%%%%%%%%%%%%%%%%%%%%%%%%%%%%%%%%%%%%%%%%%%%%%%%%%%%%%%%%%%%%%%%%%%%%%%%%%
%%%%%%%%%%%%%%%%%%%%%%%%%%%%%%%%%%%%%%%%%%%%%%%%%%%%%%%%%%%%%%%%%%%%%%%%%
%%          Reference
%%%%%%%%%%%%%%%%%%%%%%%%%%%%%%%%%%%%%%%%%%%%%%%%%%%%%%%%%%%%%%%%%%%%%%%%%%
%%%%%%%%%%%%%%%%%%%%%%%%%%%%%%%%%%%%%%%%%%%%%%%%%%%%%%%%%%%%%%%%%%%%%%%%%%\newpage


\begin{thebibliography}{99}




%\bibitem{cr}
% Cerrai, S.; R\"{o}ckner, M.; 
% Large deviations for stochastic reaction-diffusion systems with multiplicative noise and non-Lipschitz %reaction term. Ann. Probab. 32 (2004), no. 1B, 1100--1139.


%\bibitem{cm}
% Chenal, F.; Millet, A.; 
% Uniform large deviations for parabolic SPDEs and applications. 
% Stochastic Process. Appl. 72 (1997), no. 2, 161--186.

%\bibitem{ch}
% Chow, P. L.; 
% Large deviation problem for some parabolic It? equations. Comm. Pure Appl. Math. 45 (1992), 
% no. 1, 97--120.

\bibitem{dz}
 Dembo, A.; Zeitouni, O.;
 Large deviations techniques and applications. Second edition. Applications of Mathematics (New York), 38. Springer-Verlag, New York, 1998. 

\bibitem{der}
 Dereich, S.; Rough paths analysis of general Banach space-valued Wiener processes. 
 J. Funct. Anal. 258 (2010), no. 9, 2910--2936.
 
 \bibitem{dd}
Dereich, S.; Dimitroff, G.; 
A support theorem and a large deviation principle for Kunita flows. 
Stoch. Dyn. 12 (2012), no. 3, 1150022, 16 pp. 

%\bibitem{dm}
% Duan, J.; Millet, A.;
%Large deviations for the Boussinesq equations under random influences. 
%Stochastic Process. Appl. 119 (2009), no. 6, 2052--2081.  
 
% \bibitem{fw} 
%  Freidlin, M. I.; Wentzell, A. D.; 
%  Reaction-diffusion equations with randomly perturbed boundary conditions. 
%  Ann. Probab. 20 (1992), no. 2, 963--986. 
 
 
\bibitem{fggr}
Friz, P.; Gess, B.; Gulisashvili, A.; Riedel S.; 
Spatial rough path lifts of stochastic convolutions.
Preprint. 
arXiv:1211.0046


 
\bibitem{fv05}
Friz, P.; Victoir, N.;
Approximations of the Brownian rough path with applications to stochastic analysis.
Ann. Inst. H. Poincar\'e Probab. Statist. 41 (2005), no. 4, 703--724. 


\bibitem{fv07}
Friz, P.; Victoir, N.;
Large deviation principle for enhanced Gaussian processes.
 Ann. Inst. H. Poincar\'e Probab. Statist. 43 (2007), no. 6, 775--785. 

\bibitem{fvbk}
Friz, P.; Victoir, N..;
Multidimensional stochastic processes as rough paths.
Theory and applications. Cambridge Studies in Advanced Mathematics, 120. Cambridge University Press, Cambridge, 2010. 


%\bibitem{gs}
%Gao, H.; Sun, C.;
%Well-posedness and large deviations for the stochastic primitive equations in two space dimensions. 
%Commun. Math. Sci. 10 (2012), no. 2, 575--593. 


\bibitem{hai1}
 Hairer, M.; Rough stochastic PDEs. Comm. Pure Appl. Math. 64 (2011), no. 11, 1547--1585.


\bibitem{hai2}
Hairer, M.;
Solving the KPZ equation. To appear in Ann. of Math.  arXiv:1109.6811
 
 
\bibitem{hw}  
Hairer, M.;  Weber, H.;
Rough Burgers-like equations with multiplicative noise. 
To appear in Probab. Theory Related Fields.  arXiv:1012.1236
 
   
\bibitem{hmw}
Hairer, M.; Maas, J.; Weber, H.;
Approximating rough stochastic PDEs. Preprint. arXiv:1202.3094

\bibitem{ina}
Inahama, Y.; 
Large deviation principle of Freidlin-Wentzell type for pinned diffusion processes.
Preprint. arXiv:1203.5177 


\bibitem{ik}
 Inahama, Y.; Kawabi, H.;
 Large deviations for heat kernel measures on loop spaces via rough paths. 
 J. London Math. Soc. (2) 73 (2006), no. 3, 797--816. 



%\bibitem{kx}
%Kallianpur, G.; Xiong, J.;
%Large deviations for a class of stochastic partial differential equations. 
%Ann. Probab. 24 (1996), no. 1, 320--345. 

\bibitem{ku}
 Kunita, H.; 
 Stochastic flows and stochastic differential equations. 
 Cambridge Studies in Advanced Mathematics, 24. Cambridge University Press,  1990. 



\bibitem{lqz}
Ledoux, M.; Qian, Z.; Zhang, T.;
Large deviations and support theorem for diffusion processes via rough paths. 
Stochastic Process. Appl. 102 (2002), no. 2, 265--283. 

\bibitem{lcl}
 Lyons, T.; Caruana, M.; L\'evy, T.;
Differential equations driven by rough paths. 
 Lecture Notes in Mathematics, 1908. Springer, Berlin, 2007.


\bibitem{lq}
 Lyons, T.; Qian, Z.;
 System control and rough paths. Oxford Mathematical Monographs. Oxford Science Publications. 
 Oxford University Press, Oxford, 2002. 


\bibitem{ms}
 Millet, A.; Sanz-Sol\'e, M.;
 Large deviations for rough paths of the fractional Brownian motion. 
 Ann. Inst. H. Poincar\'e Probab. Statist. 42 (2006), no. 2, 245--271.


\bibitem{maas}
Maas, J.;
Malliavin calculus and decoupling inequalities in Banach spaces. 
J. Math. Anal. Appl. 363 (2010), no. 2, 383--398. 

%\bibitem{pe}
% Peszat, S.;
% Large deviation principle for stochastic evolution equations. 
% Probab. Theory Related Fields 98 (1994), no. 1, 113--136.


\bibitem{r}
Riedel, S.; Xu, W.;
A simple proof of distance bounds for Gaussian rough paths.
Preprint. arXiv:1206.5866

%\bibitem{so}
% Sowers, R. B.; 
% Large deviations for a reaction-diffusion equation with non-Gaussian perturbations. 
% Ann. Probab. 20 (1992), no. 1, 504--537. 

%\bibitem{za}
% Zabczyk, J.; 
% On large deviations for stochastic evolution equations. 
% Stochastic systems and optimization (Warsaw, 1988), 240--253, 
% Lecture Notes in Control and Inform. Sci., 136, Springer, Berlin, 1989.

%\bibitem{zh}
% Zhang, T.; 
% Large deviations for stochastic nonlinear beam equations. 
% J. Funct. Anal. 248 (2007), no. 1, 175--201.

%$\|  af\|$, 
%$\baro$
%$\interleave$

\end{thebibliography}
\end{document}